\documentclass[reqno,11pt]{amsart}
\usepackage{amscd,amssymb,verbatim}
\usepackage{hyperref}\usepackage[small,nohug,heads=vee]{diagrams}
\diagramstyle[labelstyle=\scriptstyle]

\numberwithin{equation}{section}
\theoremstyle{plain}

\newcommand\alp{\alpha}         
\newcommand\bet{\beta}
\newcommand\gam{\gamma}         \newcommand\Gam{\Gamma}
\newcommand\del{\delta}         \newcommand\Del{\Delta}
\newcommand\eps{\varepsilon}

\newcommand\tet{\theta}         \newcommand\Tet{\Theta}
\newcommand\iot{\iota}
\newcommand\kap{\kappa}
\newcommand\lam{\lambda}                \newcommand\Lam{\Lambda}

\newcommand\ome{\omega}         \newcommand\Ome{\Omega}

\newcommand\calE{{\mathcal{E}}}
\newcommand\calF{{\mathcal{F}}}

\newcommand\calH{{\mathcal{H}}}

\newcommand\calL{{\mathcal{L}}}

\newcommand\bfu{{\mathbf u}}            
\newcommand\bfv{{\mathbf v}}

\newcommand\RR{\mathbb{R}}

\newcommand\ZZ{\mathbb{Z}}

\newcommand\CC{\mathbb{C}}

\newcommand\NN{\mathbb{N}}

 \newcommand\grg{{\mathfrak{g}}}

\newcommand\nek{,\ldots,}
\newcommand\sdp{\times \hskip -0.3em {\raise 0.3ex
\hbox{$\scriptscriptstyle |$}}} 

\newcommand\End{\operatorname{End\,}}

\newcommand\Hom{\operatorname {Hom}}
\newcommand\Id{\operatorname {Id}}

\newcommand\IM{\operatorname{Im}}

\newcommand\ind{\operatorname{ind}}

\newcommand\Irr{\operatorname{Irr\, }}

\newcommand\Ker{\operatorname{Ker}}

\newcommand\ov{{\overline{v}}}

\newcommand\ow{{\overline{w}}}

\newcommand\oz{{\overline{z}}}

\newcommand\oPhi{{\overline{\Phi}}}

\newcommand\tilmu{{\widetilde{\mu}}}

\renewcommand{\>}{\rangle}
\newcommand{\<}{\langle}

\theoremstyle{plain}
\newtheorem{Thm}[subsection]{Theorem}
\newtheorem{Cor}[subsection]{Corollary}
\newtheorem{Lem}[subsection]{Lemma}
\newtheorem{Prop}[subsection]{Proposition}
\newtheorem{Conjec}[subsection]{Conjecture}

\newtheorem{Def}[subsection]{Definition}

\theoremstyle{remark}

\newtheorem{Rem}[subsection]{Remark}

\def\TeXref#1{%
        \leavevmode\vadjust{\setbox0=\hbox{{\tt
                \  {\tiny \textrm #1}}}%
        \theight=\ht0
        \advance\theight by \lineskip
        \kern -\theight \vbox to
        \theight{\rightline{\rlap{\box0}}%
        \vss}%
        }}%
\newif\ifShowLabels
\ShowLabelstrue
\newdimen\theight
\def\TeXrefEq#1{%
        \leavevmode\vadjust{\setbox0=\hbox{{\tt
                \  {\tiny \textrm #1}}}%
        \theight=\ht1
        \advance\theight by \lineskip
        \kern -\theight \vbox to
        \theight{\rightline{\rlap{\box0}}%
        \vss}%
        }}%

\newcommand{\refs}[1]{Section ~\ref{S:#1}}
\newcommand{\refss}[1]{Subsection ~\ref{SS:#1}}

\newcommand{\reft}[1]{Theorem ~\ref{T:#1}}
\newcommand{\refl}[1]{Lemma ~\ref{L:#1}}
\newcommand{\refp}[1]{Proposition ~\ref{P:#1}}

\newcommand{\refd}[1]{Definition ~\ref{D:#1}}
\newcommand{\refr}[1]{Remark ~\ref{R:#1}}
\newcommand{\refe}[1]{\eqref{E:#1}}

\newenvironment{thm}[1]%
        { \begin{Thm} \label{T:#1}  \ifShowLabels \TeXref{T:#1} \fi }%
        { \end{Thm} }

\renewcommand{\th}[1]{\begin{thm}{#1}  }
\renewcommand{\eth}{\end{thm} }

\newenvironment{lemma}[1]%
        { \begin{Lem} \label{L:#1}  \ifShowLabels \TeXref{L:#1} \fi }%
        { \end{Lem} }

\newcommand{\lem}[1]{\begin{lemma}{#1} }
\newcommand{\elem}{\end{lemma}}

\newenvironment{propos}[1]%
        { \begin{Prop} \label{P:#1}  \ifShowLabels \TeXref{P:#1} \fi }%
        { \end{Prop} }

\newcommand{\prop}[1]{\begin{propos}{#1} }
\newcommand{\eprop}{\end{propos}}

\newenvironment{corol}[1]%
        { \begin{Cor} \label{C:#1}  \ifShowLabels \TeXref{C:#1} \fi }%
        { \end{Cor} }
\newcommand{\cor}[1]{\begin{corol}{#1}  }
\newcommand{\ecor}{\end{corol}}

\newenvironment{conjec}[1]%
        { \begin{Conjec} \label{Conj:#1}  \ifShowLabels \TeXref{C:#1} \fi }%
        { \end{Conjec} }
\newcommand{\conj}[1]{\begin{conjec}{#1}  }
\newcommand{\econj}{\end{conjec}}

\newenvironment{defeni}[1]%
        { \begin{Def} \label{D:#1}  \ifShowLabels \TeXref{D:#1} \fi }%
        { \end{Def} }
\newcommand{\defe}[1]{\begin{defeni}{#1}  }
\newcommand{\edefe}{\end{defeni}}

\newenvironment{remark}[1]%
        { \begin{Rem} \label{R:#1}  \ifShowLabels \TeXref{R:#1} \fi }%
        { \end{Rem} }
\newcommand{\rem}[1]{\begin{remark}{#1}}
\newcommand{\erem}{\end{remark}}

\newcommand{\eq}[1]%
        { \ifShowLabels \TeXrefEq{E:#1} \fi
           \begin{equation} \label{E:#1} }
\newcommand{\eeq}{\end{equation}}
\newcommand{\meq}[1]%
        { \ifShowLabels \TeXrefEq{E:#1} \fi
           \begin{multline} \label{E:#1} }
\newcommand{\emeq}{\end{multline}}

\newcommand{\prf}{ \begin{proof} }
\newcommand{\eprf}{ \end{proof} }
\newcommand{\Label}[1]{\label{#1}  \ifShowLabels \TeXref{#1} \fi }

\ShowLabelsfalse

\renewcommand{\d}{\text{\( \partial\)}}
\newcommand{\p}{\bar{\d}}
\renewcommand{\b}{\bullet}

\newcommand{\n}{\nabla}

\newcommand{\E}{\calE}

\renewcommand{\L}{\calL}

\newcommand{\g}{{\Gam}}
\newcommand{\gc}{{\Gam(M,C(M))}}
\newcommand{\gme}{{\Gam(M,\E)}}
\newcommand{\ha}{^{1,0}}
\newcommand{\ah}{^{0,1}}

\newcommand{\nLC}{\n^{\text{LC}}}

\renewcommand{\v}{\mathbf{v}} 
\renewcommand{\u}{\mathbf{u}}

\renewcommand{\i}{i\, }

\newcommand{\ka}{K\"ahler }\newcommand{\kae}{K\"ahler}
\newcommand{\even}{{\operatorname{even}}}
\newcommand{\odd}{{\operatorname{odd}}}
\newcommand{\bg}{{\operatorname{bg}}}

\begin{document}

\title{Background cohomology of a non-compact K\"ahler $G$-manifold}
\author[Maxim Braverman]{Maxim Braverman$^\dag$}
\address{Department of Mathematics\\
        Northeastern University   \\
        Boston, MA 02115 \\
        USA
         }

\thanks{${}^\dag$Supported in part by the NSF grant DMS-1005888.}

\begin{abstract}
For a compact Lie group $G$ we define a regularized version of the Dolbeault cohomology of a $G$-equivariant holomorphic  vector bundles over non-compact \ka manifolds. The new cohomology is infinite-dimensional, but as a representation of $G$ it decomposes into a sum of irreducible components, each of which appears in it with finite multiplicity. Thus equivariant Betti numbers are well defined. We study various properties of the new cohomology and prove that it satisfies a Kodaira-type vanishing theorem. 
\end{abstract}
\maketitle

\section{Introduction}\label{S:introduction}

If $E$ is a holomorphic vector bundle over a compact \ka manifold, the Dobeault cohomology $H^{0,\b}(M,E)$ is finite dimensional and has a lot of nice properties. If $M$ is non-compact,  $H^{0,\b}(M,E)$ is infinite dimensional space and much less is known about it. In this paper we consider a Hamiltonian action of a compact Lie group $G$ on a non-compact manifold $M$ and assume that $E$ is a $G$-equivariant holomorphic vector bundle over $M$. If the moment map $\mu$ for this action is proper and the vector field induced by $\mu$ does not vanish outside of a compact subset of $M$, we construct a new regularized Dolbeault cohomology space $H^{0,\b}_\bg(M,E)$, called the {\em background cohomology} of $E$. It is still infinite dimensional. But as a representation of $G$ it decomposes into a direct sum of irreducible components and each component appears in this decomposition finitely many times:
\eq{finiteregIntr}
        H^{0,p}_\bg(M,E)  \ = \ \sum_{V\in \Irr G}\, \bet^p_{\bg,V}\cdot V, \qquad p=0\nek n.
\end{equation}
The alternating sum of the background cohomology is equal to the regularized index of the pair $(E,\mu)$ which was introduced in \cite{Br-index} (see also \cite{Paradan03}, \cite{MaZhang_TrIndex12} for more details and \cite{MaZhang-noncompact} for an application of the regularized index to a proof of a conjecture of Vergne \cite{Vergne07}).

The background cohomology \refe{finiteregIntr} behaves in many respects as the Dolbeault cohomology of a compact \ka  manifold. In \refs{vanishing} we prove an analogue of the Kodaira vanishing theorem for the regularized cohomology. In a separate paper \cite{BrDeformedCohCircle} we specialize to the case when $G=S^1$ is a circle group. In this case we prove an analogue of the holomorphic Morse inequalities of Witten \cite{Witten84} (see also \cite{WuZhang}).

\subsection{The assumptions}\label{INassum}
The construction of the background cohomology in this paper is done under the following two  assumptions:
\begin{enumerate}
\item The moment map $\mu$ is proper;
\item Via a $G$-invariant scalar product on the Lie algebra $\grg$ of $G$, $\mu$ induces a map $\v:M\to \grg$. Let $v$ denote the vector field on $M$ associated to this map, cf. \refe{v}. We assume that this vector field does not vanish outside of a compact subset $K$ of $M$. 
\end{enumerate}

The assumption (i) above is rather restrictive. It excludes, for example, the action of the circle group $S^1$ on $\CC^n$, which has both positive and negative weights. Unfortunately it is not clear how to define the regularized cohomology without this condition for the general compact Lie group $G$. However, in \cite{BrDeformedCohCircle} we consider the case when $G=S^1$ is a circle group and in this case extend the definition of the background cohomology to the situation when the moment map is not necessarily proper.

\subsection{The deformed cohomology}\label{SS:deformed}
A smooth strictly increasing function $s:[0,\infty)\to [0,\infty)$ is called {\em admissible} if it satisfies a rather technical growth condition at infinity, cf. \refd{Kahleradmis}.  For an admissible function $s$ we set
\[
	\phi(x) \ := \ s\big(\,|\mu(x)|^2/2\,\big)\qquad x\in M,
\]
and consider the {\em deformed Dolbeault differential}
\[
	\p_s  = \ e^{-\phi}\circ\p\circ e^\phi.
\]
We view $\p_s$ as a densely defined operator on the space $L_2\Ome^{0,p}(M,E)$ of square-integrable differential forms with values in $E$ and we define the  {\em deformed Dolbeault cohomology $H^{0,\b}_s(M,E)$} as the reduced cohomology of $\p_s$:
\[
	H^{0,p}_s(M,E) \ = \ \frac{
	\Ker\,\big(\,\p_s: L_2\Ome^{0,p}(M,E)\to L_2\Ome^{0,p+1}(M,E)\,\big)}
	{ \overline{\IM \big(\, \p_s: L_2\Ome^{0,p-1}(M,E)\to L_2\Ome^{0,p}(M,E)\,\big)}}.
\]
The space $H^{0,p}_s(M,E)$ decomposes as a sum of irreducible representations of $G$.  In \reft{finitereg} we show that 
\[
        H^{0,p}_s(M,E)  \ = \ \sum_{V\in \Irr G}\, \bet^p_{s,V}\cdot V.
\]
In other word, each irreducible representation of $G$ appears in $H^{0,p}_s(M,E)$ with finite multiplicity. 

We set 
\[
	H^{0,p}_{s,V}(M,E)\ := \  \bet^p_{s,V}\cdot V,
\]
and call it the {\em $V$-component} of the deformed cohomology.

\subsection{The background cohomology}\label{SS:Ibackground}
The multiplicities $\bet^p_{s,V}$ are non-negative integers, which depend on the choice of the admissible function $s$. The function $s$ is called {\em $V$-generic} if the value of $\bet^p_{s,V}$ is the minimal possible.  Theorems~\ref{T:independenceofs} and \ref{T:indepofh} show that for any two $V$-generic functions $s_1$ and $s_2$ the $V$-components  $H^{0,p}_{s_1,V}(M,E)$ and $H^{0,p}_{s_2,V}(M,E)$ of the deformed cohomology are naturally isomorphic.  Thus we can define the $V$-component $H^{0,p}_{\bg,V}(M,E)$ of the {\em background cohomology}  as $H^{0,p}_{s,V}(M,E)$ for some $V$-generic function $s$. 

The {\em background cohomology} of $E$ is by definition the direct sum
\[
	H^{0,p}_{\bg}(M,E)\ := \ \sum_{V\in \Irr G}\, H^{0,p}_{\bg,V}(M,E).
\]

\subsection{Kodaira-type vanishing theorem}\label{SS:Intkodaira}
Let $L$ be a positive $G$-equivariant line bundle over $M$. In \refs{vanishing} we prove the following extension of the Kodaira vanishing theorem to our non-compact setting: for every irreducible representation $V$ of $G$ there exists a integer $k_0>0$, such that for all $k\ge k_0$ the $V$-component of the background cohomology 
\[
	H^{0,p}_{\bg,V}(M,E\otimes L^{\otimes k}) \ = \ 0,
\]
for all $p>0$. 

\subsection{The organization of the paper}\label{SS:organ} This paper is organized as follows: In \refs{notation} we introduced some notations used in the rest of the paper. In \refs{index} we recall some results about the index theory on non-compact manifolds constructed in \cite{Br-index}. In \refs{tamed} we introduce the notion of an admissible function and show that the set of admissible functions is a non-empty convex cone. By technical reasons we also extend our construction to manifolds which  are only asymptotically \kae. In \refs{regularizedcoh} we define the background cohomology. We also compute the background cohomology in a simple but important example.  In Sections~\ref{S:prindependenceofs} and \ref{S:prindepofh} we prove that the background cohomology is independent of all the choices made in the definition. In \refs{vanishing} we prove a Kodaira-type vanishing theorem for the background cohomology. 

\subsection*{Acknowledgment} I am grateful to the referee for correcting several mistakes and misprints in the preliminary version of this manuscripts and for many useful remarks and suggestions. 

\section{Preliminaries and notations}\label{S:notation}

\subsection{Complex manifolds}\label{SS:complexman}
Throughout the paper $M$ denotes a complex manifold without boundary endowed with a Riemannian metric $g^M$, such that the complex structure  
\[
	J:\,TM\ \to\ TM
\] 
is an anti-symmetric operator, $J^*=-J$. The complexification  $TM\otimes\CC$ of $TM$ decomposes into the direct sum of the holomorphic and the antiholomorphic tangent bundles:
 \begin{equation*}
 	TM\otimes\CC \ = \ T^{1,0}M\oplus T^{0,1}M,
 \end{equation*}
 where $J|_{T^{1,0}M}=i, \ J|_{T^{0,1}M}=-i$. Each vector $v\in TM$ has a unique decomposition 
 \begin{equation*}
 	v \ = \ v^{1,0}\ + \ v^{0,1}, \qquad v^{1,0}\in T^{1,0}M, \ \ v^{0,1}\in T^{0,1}M.
 \end{equation*} 
Similarly, the cotangent bundle $T^*M$ decomposes as 
\begin{equation*}
 	T^*M\otimes\CC \ = \ \big(T^{1,0}M\big)^*\oplus \big(T^{0,1}M\big)^*.
 \end{equation*}

Let $I:TM\to T^*M$ denote the isomorphism defined by the Riemannian metric $g^M$. Then for any two vector fields $u$ and $v$ on $M$ we have
\eq{Iu(v)}
	Iu(v)\ =\ g^M(u,v).
\end{equation}
Notice that 
\[
	I(T^{1,0}M) \ = \ \big(T^{0,1}M\big)^*, \quad I(T^{0,1}M) \ = \ \big(T^{1,0}M\big)^*
\]
For a smooth function $f:M\to \CC$ we write
\eq{partial}
	\partial f \ = \ I\big(\,\n f^{0,1}\,\big), \qquad 
	\p f \ = \ I\big(\,\n f^{1,0}\,\big).
\end{equation} 
Then the de Rham differential $d=\partial+\p$.

\subsection{Hamiltonian group action}\label{SS:hamiltonianaction}
Suppose now the the metric $g^M$ is \ka and let $\ome$ denote the \ka form on $M$.  Our sign convention is that 
\[
	\ome(u,Jv) \ = \ g(u,v).
\]
From \refe{Iu(v)} we see that
\[
	\iot_u\ome\ = \ I(Ju), \qquad u\in TM.
\]

Assume that a compact Lie group $G$ acts on $M$ preserving the complex structure $J$ and the Riemannian metric $g^M$. Let $\grg$ denotes the Lie algebra of $G$. 

A vector $\bfu\in \grg$ generates a vector field
\eq{generatingfield}
	\bfu_M(x) \ := \  \frac{d}{dt}\big|_{t=0}\,\exp(t\bfu)\cdot x, \qquad x\in M,
\end{equation}
on $M$. The action of $G$ on $M$ is {\em Hamiltonian} if their exists a {\em moment map} $\mu:M\to \grg^*$ such that 
\eq{dmu=}
	d\<\mu,\bfu\>\ = \ \iot_{\bfu_M}\ome\ = \  I(J\bfu_M).
\end{equation}
Equivalently,
\eq{bfv=Jmu}
	\bfu_M \ =\ -\,J\,\n\<\mu,\bfu\>.
\end{equation}
Using \refe{partial} we obtain
\eq{bfe=Idmu}
	\bfu^{1,0}_M \ = \ 
	         -\,i\, I^{-1}\p\<\mu,\bfu\>, \qquad \bfu^{0,1}\ = \ i\, I^{-1}\partial \<\mu,\bfu\>. 
\end{equation}

\subsection{The Clifford action and the Dirac operator}\label{SS:clifford}
Let $\Lam^{0,\b}\ := \ \Lam^\b(T^{0,1}M)^*$ denote the bundle of antiholomorphic forms over $M$. It is endowed with a natural Clifford action of $TM\simeq T^*M$, defined by
\eq{Clifford}
	c(v)\alp \ = \ \sqrt{2}\,\big(\,Iv^{1,0}\wedge\alp\ - \ \iot_{v^{0,1}}\alp\,\big).
\end{equation}

Let $E$   be holomorphic vector bundle over $M$. We extend the Clifford action \refe{Clifford} to the tensor product 
\eq{spinorbundle}
	\E \ := \ E\otimes \Lam^\b(T^{0,1}M)^*
\end{equation}
The space of smooth section of $\E$ is denoted by $\Ome^{0,\b}(M,E)$ and is called the space of antiholomorphic differential forms on $M$ with values in $E$. The holomorphic structure on $E$ defines the antiholomorphic differential $\p:\Ome^{0,\b}(M,E)\to \Ome^{0,\b+1}(M,E)$ and, hence, the Dolbeault complex 
\[
	0\to \Ome^{0,0}(M,E) \ \overset{\p}\longrightarrow \ \Ome^{0,1}(M,E)  \overset{\p}\longrightarrow\cdots \overset{\p}\longrightarrow \Ome^{0,n}(M,E)\ \to 0.
\]
The cohomology of this complex is called the {\em Dolbeault cohomology} of $E$ and is denoted by $H^{0,\b}(M,E)$.

Fix a Hermitian metric $h^E$ and a holomorphic connection $\n^E$ on $E$. The metrics $g^M$ and $h^E$ define the $L^2$-metric on $\Ome^{0,\b}(M,E)$. Let 
\[
	\p^*:\Ome^{0,\b}(M,E)\to \Ome^{0,\b-1}(M,E)
\]
denote the adjoint operator of $\p$ with respect to this metric. We consider the {\em Dolbeault-Dirac operator} 
\eq{Dirac}
	D \ = \ \sqrt{2}\,\big(\p+\p^*\big):\,\Ome^{0,\b}(M,E)\ \to\ \Ome^{0,\b}(M,E).
\end{equation}
If $M$ is a closed manifold then the kernel of $D$ is naturally isomorphic to the Dolbeault cohomology $H^{0,\b}(M,E)$. 

If the metric $g^M$ is \kae, then $D$ coincides with the Dirac operator associated with the Clifford action \refe{Clifford} and the connection $\n^E$, cf. \cite[Prop.~3.67]{BeGeVe}. More precisely, let 
\eq{ncalE}
	\n^\E\ = \ \n^E\otimes1\ + \ 1\otimes \n^{\operatorname{LC}}:\, \Ome^{0,\b}(M,E)\ \to\ \Ome^{0,\b}(M,E)
\end{equation}
denote the connection on $\E$ defined as the tensor product of $\n^E$ and the Levi-Civita connection $ \n^{\operatorname{LC}}$ on $TM\simeq T^*M$. Then 
\[
	D \ = \ \sum_{j=1}^n\,c(e_i)\,\n^\E_{e_i},
\]
where $\{e_1\nek e_n\}$ is an orthonormal basis of $TM$.

\section{Index on non-compact manifolds}\Label{S:index}

In this section we recall the construction of the analytical index of a Dirac operator on a non-compact manifold $M$, which was introduced in \cite{Br-index} and further studied in \cite{MaZhang_TrIndex12}. It was also shown in \cite{Br-index} that this analytical index is equal to the topological index of a transversally elliptic symbol studied in \cite{Paradan03}. This index was used by X.~Ma and W.~Zhang \cite{MaZhang-noncompact} in their proof of a conjecture of M.~Vergne \cite{Vergne07}. 

\subsection{Clifford modules and Dirac operators}\Label{SS:dirac}
First, we recall the basic properties of Clifford modules and Dirac operators. When possible, we follow the notation of
\cite{BeGeVe}.

Suppose $(M,g^M)$ is a complete Riemannian manifold without boundary.  Let $C(M)$ denote the Clifford bundle of $M$ (cf. \cite[\S3.3]{BeGeVe}), i.e., a vector bundle, whose fiber at every point $x\in M$ is isomorphic to the Clifford algebra $C(T^*_xM)$ of the cotangent
space.

Suppose $\E=\E^+\oplus\E^-$ is a $\ZZ_2$-graded complex vector bundle on $M$ endowed with a graded action
\[
        (a,s) \ \mapsto \ c(a)s, \quad \mbox{where} \quad
                        a\in \gc, \ s\in \gme,
\]
of the bundle $C(M)$. We say that $\E$ is a {\em ($\ZZ_2$-graded self-adjoint) Clifford module \/} on $M$ if it is equipped with a Hermitian metric such that the operator $c(v):\E_x\to\E_x$ is skew-adjoint, for all $x\in M$ and $v\in T_x^*M$.

A {\em Clifford connection}  on $\E$ is a Hermitian connection $\n^\E$, which preserves the subbundles $\E^\pm$ and
\[
    [\n^\E_X,c(a)] \ = \ c(\nLC_X a), \quad    \mbox{for any} \quad  a\in \gc, \ X\in\g(M,TM),
\]
where $\nLC_X$ is the Levi-Civita covariant derivative on $C(M)$ associated with the Riemannian metric on $M$.

The {\em Dirac operator \/} $D:\gme\to\gme$ associated to a Clifford connection $\n^\E$ is defined by the following composition
\[
  \begin{CD}
        \gme @>\n^\E>> \g(M,T^*M\otimes \E) @>c>> \gme.
  \end{CD}
\]
In local coordinates, this operator may be written as $D=\sum\,c(dx^i)\,\n^\E_{\d_i}$. Note that $D$ sends even sections
to odd sections and vice versa: $D:\, \Gam(M,\E^\pm)\to \Gam(M,\E^\mp)$.

Consider the $L^2$-scalar product on the space of sections $\gme$ defined by the Riemannian metric on $M$ and the Hermitian structure on $\E$. By \cite[Proposition~3.44]{BeGeVe}, the Dirac operator associated to a Clifford connection $\n^\E$ is formally self-adjoint with respect to this scalar product. Moreover, it is essentially self-adjoint with the initial domain smooth, compactly supported sections, cf. \cite{Chernoff73}, \cite[Th.~1.17]{GromLaw83}.

\subsection{Group action. The index.}\Label{SS:G}
Suppose that a compact Lie group $G$ acts on $M$. Assume that there is given a lift of this action to $\E$, which preserves the grading, the connection, and the Hermitian metric on $\E$. Then the Dirac operator $D$ commutes with the action of $G$. Hence, $\Ker D$ is a $G$-invariant subspace of the space $L^2(M,\E)$ of square-integrable sections of $\E$.

If $M$ is compact, then $\Ker D^\pm$ is finite dimensional. Hence, it breaks into a finite sum \/ $\Ker D^\pm = \sum_{V\in \Irr G}\, m^\pm_V\, V$, where the sum is taken over the set $\Irr G$ of all irreducible representations of $G$. This allows one to defined the {\em index}
\eq{char}
        \ind_G(D) \ = \ \sum_{V\in \Irr G}\, (m^+_V-m^-_V)\cdot V,
\end{equation}
as a virtual representation of $G$.

Unlike the numbers $m^\pm_V$, the differences $m^+_V-m^-_V$ do not depend on the choice of the connection $\n^\E$ and the metric $h^\E$. Hence, the index $\ind_G(D)$ depends only on $M$ and the equivariant Clifford module $\E=\E^+\oplus\E^-$. We set $\ind_G(\E):= \ind_G(D)$,and refer to it as the {\em index} of $\E$.

\subsection{A tamed non-compact manifold} \Label{SS:tamed}
In \cite{Br-index} we defined and studied an analogue of \refe{char} for a $G$-equivariant Clifford module over a complete {\em non-compact \/} manifold with an extra structure. This structure is given by an equivariant map $\v:M\to\grg$, where $\grg$ denotes the Lie algebra of $G$ and $G$ acts on it by the adjoint representation. By \refe{generatingfield} $\v$ induces a vector field $v$ on $M$ defined by
\eq{v}
    v(x) \ := \bfv(x)_M\ = \ \ \frac{d}{dt}\Big|_{t=0}\, \exp{(t\v(x))}\cdot x.
\end{equation}
\defe{tamed}
Let $M$ be a complete $G$-manifold. A {\em taming map} is a $G$-equivariant map $\v:M\to\grg$, such that the vector field $v$ on $M$, defined by \refe{v}, does not vanish anywhere outside of a compact subset of $M$. If $\v$ is a taming map, we refer to the pair $(M,\v)$ as a {\em tamed $G$-manifold}.

If, in addition, $\E$ is a $G$-equivariant $\ZZ_2$-graded self-adjoint Clifford module over $M$, we refer to the pair $(\E,\v)$ as a {\em tamed Clifford module} over $M$.
\edefe
The index we are going to define depends on the (equivalence class) of $\v$.

\subsection{A rescaling of $v$} \Label{SS:rescaling}
Our definition of the index uses certain rescaling of the vector field  $v$. By this we mean the product $f(x)v(x)$, where $f:M\to[0,\infty)$ is a smooth positive function. Roughly speaking, we demand that $f(x)v(x)$ tends to infinity ``fast enough" when $x$ tends to infinity. The precise conditions we impose on $f$ are quite technical, cf. \refd{admissible}. Luckily, our index turns out to be independent of the concrete choice of $f$. It is important, however, to know that at least one admissible function exists. This is this is proven in Lemma~1.7 of \cite{Br-index}.

We need to introduce some additional notations.

For a vector $\u\in\grg$, we denote by $\L^\E_\u$ the infinitesimal action of $\u$ on $\Gam(M,\E)$ induced by the action on $G$ on $\E$. On the other side, we can consider the covariant derivative $\n_{u}^\E:\Gam(M,\E)\to\Gam(M,\E)$ along the vector field $u$ induced by $\u$. The difference between those two operators is a bundle map, which we denote by
\eq{mu}
    \mu^\E(\u) \ := \ \n^\E_{u}-\L^\E_\u \ \in \ \End \E.
\end{equation}

We will use the same notation $|\cdot|$ for  the norms on the bundles $TM, T^*M, \E$.  Let $\End(TM)$ and $\End(\E)$ denote the bundles  of endomorphisms of $TM$ and $\E$, respectively. We will denote by $\|\cdot\|$ the norms on these bundles  induced by $|\cdot|$. To simplify the notation, set
\eq{nu}
    \nu=|\v|+\|\nLC v\|+\|\mu^\E(\v)\|+|v|+1.
\end{equation}
\defe{admissible}
We say that a smooth $G$-invariant function $f:M\to[0,\infty)$ on a tamed $G$-manifold $(M,\v)$  is {\em admissible} for the triple $(\E,\v,\n^\E)$ if
\eq{limv}
        \lim_{M\ni x\to\infty}\, \frac{f^2|v|^2}{|df||v|+f\nu+1 }  \ = \ \infty.
\end{equation}
\edefe
By Lemma~1.7 of \cite{Br-index} the set of admissible functions is not empty.
\subsection{Index on non-compact manifolds}\Label{SS:noncomind}
We use the Riemannian metric on $M$, to identify the tangent and the cotangent bundles to $M$. In particular, we consider $v$ as a section of \/ $T^*M$.

Let $f$ be an admissible function. Consider the {\em deformed Dirac operator}
\eq{Dv}
        D_{fv} \ = \ D \ + \ {\i}c(fv).
\end{equation}
This is again a $G$-invariant essentially self-adjoint operator on $M$, cf. the remark on page~411 of \cite{Chernoff73}.

One of the main results of \cite{Br-index} is the following
\th{finite}
Suppose $f$ is an admissible function.  Then

1. \ The kernel of the deformed Dirac operator $D_{fv}$ decomposes, as a Hilbert space, into an infinite direct sum
\eq{finite}
        \Ker D^\pm_{fv} \ = \ \sum_{V\in \Irr G}\, m^\pm_V\cdot V.
\end{equation}
In other words, each irreducible representation of $G$ appears in $\Ker D^\pm_{fv}$ with finite multiplicity.

2. \ The differences $m^+_V-m^-_V$ $(V\in\Irr G)$ are independent of the choices of the admissible function $f$ and the $G$-invariant Clifford connection on $\E$, used in the definition of $D$.
\eth
Following \cite{Br-index}, we refer to the pair $(D,\v)$ as a {\em tamed Dirac operator}.  The above theorem allows to define the index 
\[
	\ind_G(D,\v) \ := \ \ind_G(D_{fv})
\]
using \refe{char}.

\section{Tamed asymptotically K\"ahler manifold}\label{S:tamed}

In this section we first formulate the assumptions on a manifold under which the background cohomology is constructed in the next section.  We then define the notion of an admissible function, which will be used in the construction of the regularized cohomology. The main result of this section is that the set of admissible functions is a non-empty convex cone, cf. Lemmas~ \ref{L:admissible2} and \ref{L:admissible1}.

\subsection{Asymptotically \ka manifolds}\label{SS:asKahler}
Let $(M,g^M)$ be a complex Riemannian manifold without boundary. 
\defe{asymptKa}
We say that $(M,g^M)$ is {\em asymptotically \ka} if there exists a compact subset $K\subset M$, such that the restriction of $g^M$ to $M\backslash{}K$ is a \ka metric. 

We refer to $M\backslash{}K$ as the {\em \ka part} of $M$.  
\edefe

From now on we assume that $(M,g^M)$ is asymptotically K\"ahler. 

\subsection{Asymptotically Hamiltonian group action}\label{SS:asHamiltonian}
Suppose a compact Lie group $G$ acts holomorphically on an asymptotically \ka manifold $M$.

\defe{asymptHam}
We say that a holomorphic action of $G$ on an asymptotically \ka manifold $(M,g^M)$ is {\em asymptotically Hamiltonian} if there exists a $G$-invariant compact subset $K\subset M$, such that $M\backslash{K}$ is a \ka manifold and the restriction of the action of $G$ to $M\backslash{K}$ is Hamiltonian. In other words we assume that there exists a moment map $\mu:M\backslash{}K\to \grg^*$ such that for every $\bfu\in \grg$
\[
	d\<\mu,\bfu\> \ = \ \iot(\bfu_M)\,\ome,
\]
where $\ome\in \Ome^2(M\backslash{}K)$ is the \ka form on $M\backslash{K}$ and $\bfu_M$ is the verctor field defined in \refe{generatingfield}. 
\edefe

\subsection{Tamed asymptotically \ka manifolds}\label{SS:tamedKahler}
Fix an invariant scalar product $\<\cdot,\cdot\>$ on $\grg^*$ and let $|\mu|^2$ denote the square of the norm of $\mu$ with respect to this scalar product. This scalar product defines and isomorphism $\psi:\grg^*\to \grg$ and we denote
\eq{bfv(x)}
	\bfv(x) \ := \  \psi\big(\mu(x)\big)\ \in \grg.
\end{equation}
Consider the vector field 
\eq{vfieldv}
	v(x)\ := \ -\,J\,\n\,\frac{|\mu(x)|^2}2
\end{equation}
on $M$. By \refe{bfv=Jmu}, the restriction of $v$  to  $M\backslash{K}$  is equal to the vector field $\bfv_M$ generated by $\bfv$, cf. \refe{generatingfield}.
\defe{tamedKahler}
A {\em tamed complex manifold} is a complete asymptotically \ka manifold $(M,g^M)$ together with an asymptotically Hamiltonian action of  the group $G$, such that 
\begin{enumerate}
\item the moment map $\mu$ is proper;
\item the vector field \refe{vfieldv} does not vanish outside of a compact set.
\end{enumerate}

If, in addition, the metric $g^M$ is \ka everywhere on $M$, then we refer to $(M,g^M)$ as a {\em tamed \ka manifold}.
\edefe

\subsection{Rescaling}\label{SS:Kahlerrescaling}
Let us choose a $G$-equivariant extension $\tilmu:M\to \grg^*$ of the moment map $\mu:M\backslash K\to \RR$.  As in \refe{bfv(x)}, \refe{vfieldv} we set
\eq{tilbfv(x)}
	\bfv(x) \ := \  \psi\big(\tilmu(x)\big)\ \in \grg, \quad v(x)\ := \ -\,J\,\n\,\frac{|\tilmu(x)|^2}2, \qquad x\in M.
\end{equation}

Our definition of the regularized cohomology uses certain rescaling of the function  $|\tilmu|^2/2$. By this we mean the function 
\eq{phi(x)}
	\phi(x) \ :=  \ s\big(\,|\tilmu(x)|^2/2\,\big),
\end{equation}
where $s:[0,\infty)\to [0,\infty)$ is a smooth positive strictly increasing function satisfying certain growth conditions at infinity. Roughly speaking, we demand that $\phi(x)$ tends to infinity ``fast enough" when $x$ tends to infinity. The precise conditions we impose on $\phi$ are quite technical, cf. \refd{Kahleradmis}, but our construction turns out to be independent of the concrete choice of $\phi$. It is important, however, to know that at least one admissible function exists.  This is proven in \refl{admissible2} below. 

The construction of the rescaling in this section is an adaptation of the rescaling procedure for a taming vector field in  \cite{Br-index}. In particular, the  condition we impose on the rescaling function $s$ is similar to the one in Section~2.5 of \cite{Br-index}, see also \refd{admissible} above. To formulate it we use the notation of \refss{rescaling}. More precisely we consider a $G$-equivariant holomorphic vector bundle $E$  over $M$ endowed with  a $G$-invariant holomorphic connection $\n^E$. Fix a  Hermitian metric $h^E$ on $E$. Consider the bundle $\E= E\otimes\Lam^\b(T\ah{}M)^*$ and endow it with the Hermitian metric induced by $h^E$ and $g^M$. Let $\n^\E$ denote the connection on $\E$ induced by $\n^E$ and the Levi-Civita connection on $TM$, cf. \refe{ncalE}. Then we define the function $\nu$ by \refe{nu}.

\defe{Kahleradmis}
A smooth  function $s:[0,\infty)\to [0,\infty)$ is called {\em admissible for the quadruple $(M,g^M,E,h^E)$} if $s'(r)>0$  and the function 
\eq{Kahleradmis}
	f(x)\ :=\ s'\big(|\tilmu(x)|^2/2\big)
\end{equation} 
satisfies the following condition
\eq{admissible}
	        \lim_{M\ni x\to\infty}\, \frac{f^2|v|^2}{    |df||v|+f\nu+1 }  \ = \ \infty.
\end{equation}
 
We denote by $\calF= \calF(M,g^M,E,h^E)$ the set of admissible functions for $(M,g^M,E,h^E)$
\edefe

Clearly, the above definition is independent of the choice of the extension $\tilmu$ of $\mu$. 

\rem{fhereandfthere}
The condition \refe{admissible} on $f$ is exactly the same as \refe{limv}. Notice however that while in  \refe{limv} $f$ was an arbitrary $G$-invariant function now we demand that $f$ is a function of the square of the moment map. In fact, the assumption (i) of \refd{tamedKahler} implies that $M$ is a manifold with a cylindrical end. More precisely, there exists a compact set $K\subset M$ and a diffeomorphism  
\[
	\Phi:\,M\backslash K \ \longrightarrow N\times[1,\infty), \quad \Phi(x)=(y,t), \qquad x\in M, \ y\in N, \ t\in [1,\infty)
\] 
with   $t=|\mu|^2/2$. Thus we now require that the restriction of $f$ to the cylindrical end depends only on the coordinate $t$. 
\erem

\rem{h1<h2}
Suppose $h_1^E\le h_2^E$ are two Hermitian metrics on $E$. Then it follows immediately from the definition that $\calF(M,g^M,E,h^E_2)\subseteq \calF(M,g^M,E,h^E_1)$.
\erem
\lem{admissible2}
Given  a holomprphic Hermitian $G$-equivariant vector bundle $(E,h^E)$ over a tamed  asymptotically \ka manifold $(M,g^M)$,  the set $\calF= \calF(M,g^M,E,h^E)$ of admissible functions is not empty. Moreover, for any function $\kap:[0,\infty)\to [0,\infty)$ there exists an admissible function $s\in \calF$ such that $s(t)\ge\kap(t)$  and $s'(t)\ge\kap(t)$ for all $t>0$. 
\elem
\prf
Choose a smooth function $r:[0,\infty)\to [0,\infty)$ such that
\begin{enumerate}
\item $r(t)\ge \kap(t)$ for all $t\ge 0$ and $\lim_{t\to\infty}r(t)=\infty$;
\item $r'(t)>0$ and $\lim_{t\to\infty}r'(t)=\infty$;
\item $r''(t)>0$ for all $t\ge0$. 
\end{enumerate}

Since $\tilmu$ is proper the set $\{x\in M:\,|\tilmu(x)|=t\}$ is compact. Set
\eq{ab}
  \begin{aligned}
	a(t) \ &= \ \min\big\{\,|v(x)|^2:\, |\tilmu(x)|=t\,\big\}; \\	
	b(t) \ &= \ \max\big\{\, |d\tilmu(x)|\,|v(x)|+\nu(x)+1:\, |\tilmu(x)|=t\,\big\}.
  \end{aligned}
\end{equation}
Since $v(x)\not=0$ outside of a compact set, we conclude that $a(t)\not=0$ for large $t$. Hence, we can choose a smooth strictly increasing function $c:[0,\infty)\to [0,\infty)$ which increases fast enough so that 
\begin{enumerate}
\item $c(t)\ge\kap(t)$ and $\lim_{t\to\infty} c(t)=\infty$; \hfill \refstepcounter{equation}(\theequation)\label{E:c>kap}
\item the function $c(t)e^{-t}>1$ and is strictly increasing;
\item \(\displaystyle \lim_{t\to\infty}\,c(t)\,\frac{a(t)}{b(t)} \ = \ \infty.\) \hfill \refstepcounter{equation}(\theequation)\label{E:abc}
\item \(\displaystyle c(t)\ge \frac{r'(t)^2}{r''(t)} \) for all $t\ge 0$. \hfill \refstepcounter{equation}(\theequation)\label{E:c>r}
\end{enumerate}
Then from (ii) we conclude that for $t<\tau$ 
\[
	\frac1{c(\tau)}\ < \ \frac{e^{t-\tau}}{c(t)}.
\]
Hence,
\eq{intc}
	\int_t^\infty\,\frac{d\tau}{c(\tau)} \ \le \ \frac1{c(t)}\cdot \int_t^\infty\,e^{t-\tau}\,d\tau \ = \ \frac{1}{c(t)}.
\end{equation}
Set
\[
	s(t) \ = \ r(0)\ + \ 
	  \int_0^{t}\, \frac1{ \int_u^\infty\,\frac{d\tau}{c(\tau)} }\,du.
\]
Then   $s'(t)>0$ for all $t>0$. We will show that $s(t)\in \calF$. 

First, using \refe{intc} we obtain 
\eq{s'(t)gec(t)}
	s'(t) \ = \ \frac1{\int_{t}^\infty\,\frac{d\tau}{c(\tau)}} \ \ge \ c(t) \ \longrightarrow \ \infty, \qquad
	\text{as}\quad t\to\infty.
\end{equation}
In particular, $s'(t) \ge c(t)\ge 1$ for large $t$. It follows that 
\eq{s'>c}
	|s''(t)| \ = \ \frac{1/c(t)}{\left(\,\int_t^\infty\,\frac{d\tau}{c(\tau})\,\right)^2} \ = \ 
	\frac{s'(t)^2}{ c(t) } \ \ge \  s'(t) \ \ge\ 1.
\end{equation}
 Hence, 
 \eq{s''+s'+1}
 	|s''(t)|\ +\ s'(t)\ +\ 1 \ \le \ 3\,|s''(t)|\ = \ \frac{3\,s'(t)^2}{c(t)}.
 \end{equation}

Using \refe{s''+s'+1}, and the definition \refe{ab} of the functions $a(t)$ and $b(t)$ we conclude that for all $x\in M$ with $\tilmu(x)= t$ \ ($|t|\gg1$) we have
\[
  \begin{aligned}
	\frac{f^2|v|^2}  {    |df|\, |v|+|f|\,\nu+1 } \ &= \  
			\frac{s'(\tilmu)^2|v|^2}  {    | s''(\tilmu) |\, |d\tilmu|\,|v|+|s'(\tilmu)|\,\nu+1 } \
	\ge \ \frac{a(t)}{b(t)}\cdot \frac{s'(t)^2}{|s''(t)|+|s'(t)|+1}
	\\ &\ge \ \frac{a(t)}{b(t)}\cdot \frac{s'(t)^2}{3s'(t)^2/c(t)}
	\ = \ \frac1{3}\,  \frac{a(t)}{b(t)}\,c(t).
  \end{aligned}
\]
By \refe{abc} the right hand side of the above inequality tends to infinity as $t\to\infty$. Hence, $s(t)$  is an admissible function. 

From \refe{c>kap} and \refe{s'(t)gec(t)} we conclude that $s'(t)\ge\kap(t)$ for all $t>0$. It remains to show that $s(t)\ge\kap(t)$. From \refe{c>r} we obtain 
\[
	\frac{1}{c(t)}\ \le \ \frac{r''(t)}{r'(t)^2} \ = \ -\,\frac{d}{dt}\,\frac1{r'(t)}.
\]
Hence, 
\[
	\int^\infty_t\,\frac1{c(\tau)}\,d\tau \ \le \ \frac1{r'(t)}
\]
and 
\[
	s(t) \ = \ r(0) \ + \ 
	  \int_0^{t}\, \frac1{ \int_u^\infty\,\frac{d\tau}{c(\tau)} }\,du \ \ge \ r(t) \ \ge \kap(t).
\]
\eprf
\lem{admissible1}
Suppose $s_1$ and $s_2$ are admissible functions for the quadruple $(M,g^M,E,h^E)$.  Then for any  positive real numbers  $t_1,t_2>0$ the function $s:=t_1s_1+t_2s_2$ is admissible. Thus, the set $\calF= \calF(M,g^M,E,h^E)$ of admissible functions is a convex cone. 
\elem
\prf
Let $s_1, s_2$ be admissible functions and let $s=t_1s_1+t_2s_2$. Set 
\[
	f(x)\ := \ s'\big(\,|\tilmu(x)|^2/2\,\big), \quad f_i(x)\ := \ s'_i\big(\,|\tilmu(x)|^2/2\,\big), \ \  i=1,2.
\]
Then $f=t_1f_1+t_2f_2$. Since $f_1,f_2>0$ we have
\eq{f2>}
	f^2\ \ge \  t_1^2f_1^2+t_2^2f_2^2 \ \ge \ \min\big\{t_1^2,\, t_2^2\,\big\}\cdot\big(\,f_1^2+f_2^2\,\big).
\end{equation}
Also
\eq{denominator<}
 \begin{aligned}
	 |df||v|+|f|\nu+1 \  & \le \ t_1\,\big(\, |df_1||v|+|f_1|\nu+1\,\big)\ + \ t_2\,\big(\, |df_2||v|+|f_2|\nu+1\,\big)\\
	 \ &\le \ (t_1+t_2)\cdot
	 	\max\,\big\{\, |df_1||v|+|f_1|\nu+1,\,|df_2||v|+|f_2|\nu+1\,\big\}.
 \end{aligned}
\end{equation}
Set 
\[
	\tau \ := \ \frac{ \min\big\{t_1^2,\, t_2^2\,\big\} }{ t_1+t_2 }\ > \ 0.
\]
Then from \refe{f2>} and \refe{denominator<} we obtain
\begin{multline}\label{E:fraction>}\notag
	\frac{f^2|v|^2}{    |df||v|+|f|\nu+1 }  \ \ge \ 
	     \frac{\tau|v|^2\,(f_1^2+f_2^2)}{  \max\big\{\,  |df_i||v|+|f_i|\nu+1:\, i=1,2\,\big\} } \\
	= \ \min\left\{\,\frac{\tau|v|^2\,(f_1^2+f_2^2)}{  |df_i||v|+|f_i|\nu+1 }:\, i=1,2\,\right\}\\
	\ge \ \tau\,\min\left\{\,\frac{|v|^2\,f_i^2}{  |df_i||v|+|f_i|\nu+1 }:\, i=1,2\,\right\}
\end{multline}
Since 
\[
	 \lim_{M\ni x\to\infty}\, \frac{f_i^2|v|^2}{    |df_i||v|+|f_i|\nu+1 }  \ = \ \infty,
\]
for $i=1,2$, we conclude from the last inequality that $f$ satisfies \refe{admissible}.
\eprf

\section{The background cohomology of a tamed K\"ahler manifold}\label{S:regularizedcoh}

From now on we assume that $(M,g^M)$ is a tamed asymptotically \ka manifold endowed with an asymptotically Hamiltonian action of a compact Lie group $G$ and that $E$ is a  $G$-equivariant holomorphic vector bundle over $M$. The purpose of this section is to define the backgraund Dolbeault cohomology $H^{0,\b}_\bg(M,E)$. First we fix a Hermitian metric $h^E$ on $E$ and an admissible function $s\in \calF(M,g^M,E,h^E)$, cf. \refd{Kahleradmis}, and use it to define a deformation $\p_s$ of the Dolbeault differential, cf. \refe{deformedDobeault}. We then define the deformed cohomology $H^{0,\b}_s(M,E)$ as the reduced cohomology of $\p_s$. We use the results of \cite{Br-index}  to prove that every irreducible representation of $G$ appears in $H^{0,\b}_s(M,E)$  with finite multiplicity. In other words,
\[
	H^{0,\b}_s(M,E) \ = \ \sum_{V\in \Irr G} H^{0,\b}_{s,V}(M,E)
\]
where each $H^{0,\b}_{s,V}(M,E)$ is a finite dimensional representation of $G$, which decomposes as a direct sum of a finite number of copies of $V$. 

The function $s$ is called {\em $V$-generic} if the dimension of $H^{0,\b}_{s,V}(M,E)$ is minimal possible.  Theorems~\ref{T:independenceofs} and \ref{T:indepofh} state that $H^{0,\b}_{s,V}(M,E)$ is independent of the choice of the generic function $s$ and the Hermitian metric $h^E$. We define the background cohomology $H^{0,\b}_\bg(M,E)$ as the deformed cohomology $H^{0,\b}_s(M,E)$ for some choice of the metric $h^E$ and the generic function $s$. We conclude this section with a computation of the background cohomology of $\CC^n$ endowed with the action of the circle group $G=S^1$. 

\subsection{Deformed Dolbeault differential}\label{SS:Dolbeault}
Let $s:[0,\infty)\to [0,\infty)$ be an admissible function, cf. \refd{Kahleradmis}, and set
\[
	\phi(x) \ := \ s\big(\,|\tilmu(x)|^2/2\,\big), \qquad x\in M.
\]
Consider the deformed Dolbeault differential
\eq{deformedDobeault}
	\p_s \alp\ = \ e^{-\phi}\circ\p\circ e^\phi\,\alp 
	\ = \ \p\alp\ + \ f\,\p\big(\,|\tilmu|^2/2\,\big)\wedge\alp,
\end{equation}
where, as in \refd{Kahleradmis}, $f= s'(|\tilmu|^2/2)$.

Let $v$ be as in \refe{tilbfv(x)}.  Then from \refe{partial} we conclude that for $x\in M\backslash{K}$
\[
	\p\big(\,|\mu|^2/2\,\big) \ = \ i\,Iv^{1,0}.
\]
Hence
\[
	\p_s\,\alp \ := \ \p\,\alp\ + \ i f\,Iv\ha\,\wedge\alp, \quad 
	\p_s^*\,\alp \ = \ \p^*\,\alp \ - \   if\,\iot(v\ah)\,\alp.
\]
The deformed Dolbeault-Dirac operator is defined by 
\eq{Ds}
	D_{s} \ = \ \sqrt2\,\big(\,\p_s \, + \, \p_s^*\,\big).
\end{equation}
Comparing with \refe{Clifford} and \refe{Dirac}  we now conclude that the restriction of $D_s$ to $M\backslash{K}$ coincides with the deformed Dirac operator $D_{fv}$ defined in \refe{Dv}. Hence, 
\eq{Dfv=pp*}
	D_s \ =\ D_{fv} \ + \ R,
\end{equation}
where $R$ is a zero-oder operator supported on $K$.

\subsection{Deformed Dolbeault cohomology}\label{SS:Dolbeaultcohomology}
Let $L_2\Ome^{0,p}(M,E)$ denote the space of square-integ\-rable differential $(0,p)$-forms on $M$ with values in $E$.%
\footnote{The space  $L_2\Ome^{0,p}(M,E)$ depends on the metrics $g^M$ and $h^E$. But we omit them from the notation for simplicity.}
We view the operators $\p_s,\p_s^*,D_s$ as densely defined operators on the space  $L_2\Ome^{0,\b}(M,E)$ and we define the {\em deformed Dolbeault cohomology $H^{0,\b}_s(M,E)$} of the triple $(M,E,s)$ as the reduced cohomology of the deformed differential $\p_s$. Thus we set
\eq{deformedDolbeault}
	H^{0,p}_s(M,E) \ = \ \frac{
	\Ker\,\big(\,\p_s: L_2\Ome^{0,p}(M,E)\to L_2\Ome^{0,p+1}(M,E)\,\big)}
	{ \overline{\IM \big(\, \p_s: L_2\Ome^{0,p-1}(M,E)\to L_2\Ome^{0,p}(M,E)\,\big)}}.
\end{equation}
The space $H^{0,\b}_s(M,E)$ is naturally isomorphic to the kernel of the deformed Dirac operator $D_s$, cf., for example, \cite[(3.1.22)]{MaMarinescu_book}. From \reft{finite} and the equation \refe{Dfv=pp*} we now obtain the following
\th{Cohfinite}
Suppose $s\in \calF(M,g^M,E,h^E)$ is an admissible function.  Then the deformed Dolbeault cohomology $H^{0,p}_s(M,E)$ decomposes, as a Hilbert space, into an infinite direct sum
\eq{finitecoh}
        H^{0,p}_s(M,E)  \ = \ \sum_{V\in \Irr G}\, \bet^p_{s,V}\cdot V.
\end{equation}
In other words, each irreducible representation of $G$ appears in $H^{0,p}_s(M,E)$ with finite multiplicity.
\eth

For each irreducible representation $V$ the finite dimensional representation 
\[
	H^{0,p}_{s,V}(M,E)  \ = \ \bet^p_{s,V}\cdot V
\]
is called the {\em $V$-component} of the deformed cohomology. Then 
\[
        H^{0,p}_s(M,E)  \ = \ \sum_{V\in \Irr G}\, H^{0,p}_{s,V}(M,E).
\]
	
From Theorems~\ref{T:finite} and  \ref{T:Cohfinite}  and the equation \refe{Dfv=pp*} we now obtain the following

\prop{finitereg}
Define $\E^+= E\otimes\Lam^\even(T\ah{}M)^*, \ \E^-=E\otimes\Lam^\odd(T\ah{}M)^*$. Then $\E=\E^+\oplus\E^-$ is a graded vector bundle over $M$. Fix an admissible function $s$ and set  $f=s'\big(|\tilmu|^2/2\big)$. Let $D_{fv}:\Gam(M,\E)\to \Gam(M,\E)$ be the Dirac operator \refe{Dv} and let $m_V^\pm$ be as in \refe{finite} . Then, for every $V\in \Irr G$,
\eq{mpm-mp1}
	m_V^+-m_V^- \ = \ \sum_{j=1}^n\,(-1)^p\,\bet_{s,V}^p.
\end{equation}
\eprop

The numbers $\bet^p_{s,V}$ are non-negative integers, which depend on the choice of the admissible function $s$. 

\defe{background}
The minimal possible value of $\bet^p_{s,V}$  is called the {\em background Betti number} and is denoted by $\bet^p_{\bg,V}$:
\begin{equation}\label{E:backgroundBetti}
	\bet^p_{\bg,V}\ := \ \min\big\{\, \bet^p_{s',V}:\, s'\in \calF\,\big\}
\end{equation}
\edef
\defe{V-generic}
An admissible  function $s$ is called {\em $V$-generic} if $\bet^p_{s,V}\ = \  \bet^p_{\bg,V}$ for all $p=0\nek n$.
\edefe

\subsection{Independence of the $V$-generic function}\label{SS:indepofsh}
The main results of this paper are the following two theorems which show that the deformed cohomology is essentially independent of the choice of a generic function $s\in \calF$ and the Hermitian metric $h^E$ on $E$.

\th{independenceofs}
Let $V\in \Irr G$ be an irreducible representation of $G$. For any two $V$-generic admissible functions $s_1,s_2\in \calF(M,g^M,E,h^E)$ there exists a canonical isomorphism 
\eq{Phis1s2}
	\Phi_{s_1s_2}^V:\, H^{0,\b}_{s_1,V}(M,E) \ \longrightarrow \ H^{0,\b}_{s_2,V}(M,E),
\end{equation}
satisfying the cocycle condition
\eq{cocyclePhi}
	\Phi_{s_2s_3}^V\circ\Phi_{s_1s_2}^V \ = \ \Phi_{s_1s_3}^V.
\end{equation}

If $\tau =s_{2}-s_{1}\ge0$ is an admissible function then the isomorphism $\Phi_{s_{1}s_{2}}^V$ is induced by the map
\eq{alptoealp}
	\alp \ \mapsto \ e^{-\tau(|\tilmu|^2/2)}\,\alp,\, \qquad \ \alp\in L_2\Ome^{0,\b}(M,E).
\end{equation}
\eth
\noindent The proof is given in \refs{prindependenceofs}.

\subsection{Independence of the Hermitian metric}\label{SS:indepofh}
Suppose $h_1^E$ and $h_2^E$ are two Hermitian metrics on $E$. By \refr{h1<h2}
\[
	\calF(M,g^M,E,h_1^E)\,\cap\,\calF(M,g^M,E,h_2^E) \ \supset \calF(M,g^M,E,h_1^E+h_2^E) \ \not= \ \emptyset.
\]
Thus there exists a function $s:[0,\infty)\to [0,\infty)$ which is admissible for both Hermitian metrics $h_1^E$ and $h_2^E$. 
\th{indepofh}
Suppose $h_i^E$ $(i=1,2)$ are two Hermitian metrics on $E$. We denote by $E_i$ the vector bundle $E$ endowed with the metric $h^E_i$. Let $V\in \Irr{}G$ and let $s$ be a $V$-generic admissible function for both bundles $E_1$ and $E_2$. Then there is a natural isomorphism 
\eq{Tets1s2}
	\Tet_{h_1^Eh_2^E}^V:\, H^{0,p}_{s,V}(M,E_1) \ \longrightarrow \ H^{0,p}_{s,V}(M,E_2),
\end{equation}
satisfying the cocycle condition
\eq{cocycleTet}
	\Tet_{h_2^Eh_3^E}^V\circ\Tet_{h_1^Eh_2^E}^V \ = \ \Tet_{h_1^Eh_3^E}^V.
\end{equation}
\eth
\noindent The proof is given in \refs{prindepofh}.

Theorems~\ref{T:independenceofs} and \ref{T:indepofh} justify the following 
\defe{backgroundcohomology}
If\/ $V\in \Irr{}G$, the {\em $V$-component of the background  cohomology} $H^{0,p}_{\bg,V}(M,E)$  of the triple  $(M,g^M,E) $ is defined to be the deformed cohomology $H^{0,\b}_{s,V}(M,E)$ for any Hermitian metric $h^E$ and any $V$-generic function $s$. 

The {\em background cohomology} $H^{0,p}_\bg(M,E)$ is the direct sum
\begin{equation}
	H^{0,p}_\bg(M,E) \ := \ \sum_{V\in \Irr G} H^{0,p}_{\bg,V}(M,E).
\end{equation}
 \edefe

\th{finitereg}
Define $\E^+= E\otimes\Lam^\even(T\ah{}M)^*, \ \E^-=E\otimes\Lam^\odd(T\ah{}M)^*$. Then $\E=\E^+\oplus\E^-$ is a graded vector bundle over $M$. Fix an admissible function $s$ and set  $f=s'\big(|\tilmu|^2/2\big)$. Let $D_{fv}:\Gam(M,\E)\to \Gam(M,\E)$ be the Dirac operator \refe{Dv} and let $m_V^\pm$ be as in \refe{finite} . Then, for every $V\in \Irr G$,
\eq{mpm-mp}
	m_V^+-m_V^- \ = \ \sum_{j=1}^n\,(-1)^p\,\bet_{\bg,V}^p.
\end{equation}
\eth

\prf
From \refp{finitereg} we conclude that \refe{mpm-mp} holds if $s$ is a $V$-generic function. But from \reft{finite}.2 it follows that the left hand side of \refe{mpm-mp}  is independent of the choice of an admissible function $s$. 
\eprf

\subsection{An examples}\label{S:example}
We finish this section by computing the background cohomology in the following simple case. Suppose $G=S^1$ is the circle group and $M=\CC^n$ endowed with the standard metric. We assume that $G=S^1$ acts on $\CC^n$ by 
\begin{equation}\label{E:acionCn}
     e^{\i\tet}:(z_1\nek z_n) 
     \ \mapsto \ (e^{-\i\lam_1\tet}z_1\nek e^{-\i\lam_n\tet}z_n),
\end{equation}
where $\lam_1\nek \lam_n$ are positive integers. The moment map is given by the formula
\[
    \mu(z) \ = \ \sum_{i=1}^n\, \lam_i|z_i|^2/2 \ \in \ \RR\simeq\grg^*.
\]

 For an integer $k\in \ZZ$ we denote by $V_k$ the one dimensional representation of $G=S^1$ on which the element $e^{i\tet}\in G$ acts by
\[
	e^{i\tet}:\,u\ \mapsto \ e^{-ik\tet}v, \qquad u\in V_k\simeq\CC.
\]

Consider the line bundle $E_k=V_k\times M$ over $M=\CC^n$. It is a $G$-equivariant line bundle with the action of $G$ given by 
\[
	e^{i\tet}\cdot(u,m) \ \mapsto \ (e^{-ik\tet}u,e^{i\tet}\cdot m).
\]

\prop{S1Cn}
If \/ $G=S^1$ acts on $\CC^n$ by \refe{acionCn}, then 
\[
	H^{0,p}_\bg(\CC^n,E_k) \ = \ 0, \qquad \text{for all} \quad p>0,
\]
and   $H^{0,0}_{\bg,V_m}(\CC^n,E_k)$  ($m\in \ZZ$) is isomorphic to the space of  polynomials in $z_1\nek z_n$ spanned by the monomials $z_1^{m_1}\cdots{}z_n^{m_n}$ such that 
\[
	\sum_{i=1}^n\,\lam_i\,m_i \ = \ m-k.
\]
\eprop
\prf
Let $s:[0,\infty)\to [0,\infty)$ be a smooth strictly increasing function such that for $t\ge1$ we have $s(t)=\sqrt{2t}$. A straightforward computation shows  that $s(t)$  is an admissible function for $E_k$. Thus the deformed cohomology $H^{0,\b}_s(\CC^n,E_k)$ is isomorphic to the kernel of the deformed Dirac operator $D_s=\sqrt2(\p_s+\p_s^*)$, where 
\[
	\p_s\ = \ e^{-\mu} \circ \p\circ e^{\mu}.
\]
This kernel is computed in  \cite{Witten84} and \cite[Prop.~3.2]{WuZhang}. It follows from this computation that $H^{0,p}_s(\CC^n,E_k)$  is equal to 0 for $p>0$ and is isomorphic to 
the space of  polynomials in $z_1\nek z_n$ spanned by the monomials $z_1^{m_1}\cdots{}z_n^{m_n}$ such that 
\[
	\sum_{i=1}^n\,\lam_i\,m_i \ = \ m-k
\]
Since $H^{0,p}_s(\CC^n,E_k)$ vanishes for all but one value of $p$, it follows from \reft{finitereg} that $s$ is $V$-regular for all $V$ and, hence, $H^{0,\b}_{\bg}(M,E)= H^{0,\b}_s(M,E)$.
\eprf

\rem{Bargman}
Recall that the Bargmann space $L^2_{\text{hol}}(\CC^{n},\nu)$ is the
Hilbert space of holomorphic functions on $\CC^n$, which are  square-integrable with respect to the measure $\nu=\i^ne^{-|z|^2/2}dzd\oz$. The Bargmann space is usually considered as a quantization of $\CC^n$, cf., for example, \cite{Getzler94}, \cite[Ch.~6]{GGK-book}. The monomials
$z_1^{m_1}{\cdots}z_n^{m_n}$ form an orthogonal basis for
$L^2_{\text{hol}}(\CC^n,\nu)$, so that this space is the completion of the background cohomology space $H^{0,0}_\bg(M,E_k)$.
\erem
\rem{bundleS1}
With a little more work and using the computations of \cite{WuZhang} one obtain an analogue of \refp{S1Cn} for the case when $M$ is the total space of a holomorphic vector bundle over a \ka manifold $B$ on which $G=S^1$ acts by fiberwise linear transformations. 
\erem
\rem{furure}
It is natural and very important to consider a more general action of $G=S^1$ on $\CC^n$ given by 
\begin{equation}\label{E:acionCn2}
     e^{\i\tet}:(z_1\nek z_r, z_{r+1}\nek z_n) 
     \ \mapsto \ (e^{-\i\lam_1\tet}z_1\nek e^{-\i\lam_r\tet}z_r,
          e^{\i\lam_{r+1}\tet}z_{r+1}\nek e^{\i\lam_n\tet}z_n),
\end{equation}
where $\lam_1\nek \lam_n$ are positive integers. Unfortunately, the moment map for this action is not proper and we can not construct the background cohomology for such an action using the methods of this paper.  In \cite{BrDeformedCohCircle} we study the background cohomology in the case $G=S^1$ and show that in this case it  can be defined even if the moment map is not proper. In particular, we compute the background cohomology for the action \refe{acionCn2}.
\erem

\section{Proof of \reft{independenceofs}}\label{S:prindependenceofs}

In this section we construct the canonical isomorphism $\Phi_{s_1s_2}^V$, cf. \refe{Phis1s2}. 
\subsection{A family of admissible function}\label{SS:backgroundharmonic}
Let $s_1,s_2\in \calF$ be admissible $V$-generic functions. Set 
\begin{equation}\label{E:s(t)}
	s(t)\ := \ s_1\ + \ t s_2.
\end{equation}
By Lemma~\ref{L:admissible1}, $s(t)\in \calF$ for all $t\ge0$.  Consider the deformed Dirac operator
\begin{equation}\label{D(t)}
		D(t) \ = \ \sqrt2\,\big(\,\p_{s(t)}\,+\,\p_{s(t)}^*\,\big).
\end{equation}
Let 
\[
	L_2\Ome^{0,\b}(M,E)^V\ := \ \Hom_G\big(\,L_2\Ome^{0,\b}(M,E),V\,\big)\otimes V
\]
denote the $V$-component of the space of square-integrable differential forms. Let $D^V(t)$ denote the restriction of $D(t)$ to $L_2\Ome^{0,\b}(M,E)^V$.

Set
\begin{equation}\label{E:Ht}
	\calH^{0,p}_t\ := \ \Ker D(t)\,\cap L_2\Ome^{0,p}(M,E)^V, \qquad p=0\nek n.
\end{equation}
Then $\calH^{0,p}_t$ is canonically isomorphic to $H^{0,p}_{s(t),V}(M,E)$. 

The following proposition states that though $D^V(t)$ is a family of unbounded operators, its kernel behave like one would expect from the kernel of an analytic family of bounded operators. 
\prop{family of projections}
There exists a discrete sequence of positive numbers $t_1, t_2, ...$ such that 
\begin{equation}\label{E:generic}
		\dim\calH^{0,p}_t\ = \ \bet^p_{\bg,V}, \qquad p=0\nek n,
\end{equation}
for all $t\ge0,\ t\not\in \{t_1,t_2,...\}$. 

Furthermore, there exists a unique continuous family of orthogonal projections 
\[
	P_t:\,L_2\Ome^{0,\b}(M,E)^V\ \to \ L_2\Ome^{0,\b}(M,E)^V, \qquad t\ge0
\]
such that 
\begin{equation}\label{E:family of projections}
	\IM P_t\ = \ \calH^{0,\b}_t,
\end{equation}
for all $t\not\in\{t_1,t_2,...\}$.
\eprop
The proof, based on the Kato's theory of holomorphic family of operators \cite[Ch.~VII]{Kato}, is given in Appendix~\ref{S:holomorphicB}. 

 Set 
\begin{equation}\label{E:barcalH(t)}
	\bar\calH^{0,\b}_t \ := \ \IM P_t, \qquad t\ge0.
\end{equation}
Then 
\[
	\bar\calH^{0,\b}_t\ \subset\ \calH^{0,\b}_t \ \simeq \ H_{s(t),V}^{0,\b}(M,E)
\] 

\subsection{A map from $\bar\calH^{0,\b}_{s(t_1)}$  to $\bar\calH^{0,\b}_{s(t_2)}$ for $t_2>t_1$}
For $t\ge0$ consider the map
\[
	\Psi_{t}:\,L_2\Ome^{0,\b}(M,E)\ \longrightarrow \ L_2\Ome^{0,\b}(M,E), \qquad
	\Psi_{t}(\alp)\ := \ e^{-ts_2(|\tilmu|^2/2)}\alp.
\]
For every $t_2>t_1\ge0$ we have
\[
	\p_{s(t_2)}\cdot\Psi_{t_2-t_2} \ = \ \Psi_{t_2-t_1}\cdot \p_{s(t_1)}.
\]
Hence, $\Psi_{t_2-t_1}$ induces a map of cohomology 
\[
	\Phi_{t_1,t_2}:\,H^{0,\b}_{s(t_1),V}(M,E) \to H^{0,\b}_{s(t_2),V}(M,E).
\]
Clearly, if $t_3>t_2>t_1\ge0$, then
\begin{equation}\label{E:Phit1t2t3}
	\Phi_{t_2,t_3}\circ \Phi_{t_1,t_2}\ = \ \Phi_{t_1,t_3}.
\end{equation}

Consider the map
\eq{Phit}
	\oPhi_{t_1,t_2}:\, \bar\calH_{s(t_1)}^{0,\b}\ \longrightarrow\ \bar\calH_{s(t_2)}^{0,\b},  \qquad 
	\oPhi_{t_1,t_2}\ := \ P_{t_2}\circ \Psi_{t_2-t_1}.
\end{equation}
Note that if $s(t_2)$ is $V$-generic and $h\in \bar\calH_{s(t_1)}^{0,\b}$ , then $\oPhi_{t_1,t_2}(h)$ is the harmonic representative of the cohomology class of $h$. Hence, it follows from \eqref{E:Phit1t2t3} that if $s(t_2)$ and $s(t_3)$ are $V$-generic then 
\begin{equation}\label{E:oPhit1t2t3}
		\oPhi_{t_2,t_3}\circ \oPhi_{t_1,t_2}\ = \ \oPhi_{t_1,t_3}.
\end{equation}

\prop{Phiinjective}
For every $t_2>t_1\ge0$ the map \refe{Phit} is bijective. 
\eprop
We present a proof of the proposition in \refss{prPhiinjective} after some additional constructions are introduced. 
\subsection{A covariant derivative on $\bar\calH^{0,p}_t$}\label{SS:nt}
We consider the collection of spaces $\bar\calH^{0,\b}_t$ $(t\ge0)$ as a finite dimensional vector bundle over $\RR_{\ge0}$. Let us define a connection on this bundle by 
\begin{equation}\label{E:nt}
 	\n_{d/dt}h(t)\ := \ P_t\,\big(\,h'(t)+s_2h(t)\,\big).
\end{equation}

\begin{Lem}\label{L:h(t)}
Fix $h_0\in \bar\calH^{0,\b}_0$ and set 
\begin{equation}\label{E:h(t)}
	h(t)\ := \ \oPhi_{0,t}\big(h_0\big)\ = \ P_t\,\big(e^{-ts_2(|\tilmu|^2/2)}h_0\big).
\end{equation}
Then $h(t)$ is a flat section, i.e., 
\begin{equation}\label{E:nh(t)=0}
	\n_{d/dt}h(t)\ = \ 0.
\end{equation}
\end{Lem}
\prf
It suffices to check \refe{nh(t)=0} for generic $t$. Fix  a $V$-generic $t_0$. Then there exists an $\eps>0$ such that all $t\in (t_0-\eps,t_0+\eps)$ are $V$-generic.  It follows from \refe{oPhit1t2t3} that for $t\in (t_0,t_0+\eps)$
\[
	\oPhi_{0,t}\ = \ \oPhi_{t_0,t}\circ\oPhi_{0,t_0}.
\] 
Thus there exists a function $\alp:(t_0,t_0+\eps)\to L_2\Ome^{0,\b}(M,E)^V$ such that $\alp(t_0)=0$ and 
\[
	h(t)\ = \ P_t\,\big(\,e^{-(t-t_0)s_2(|\tilmu|^2/2)}h(t_0)\,) \ = \ e^{-(t-t_0)s_2(|\tilmu|^2/2)}h(t_0) + \p_t\alp(t)
\]
Hence,
\[
	h'(t_0)\ = \ -s_2h(t_0)+\big(\frac{d}{dt}\big|_{t=t_0}\p_t\big)\,\alp(t_0) +   \p_t\alp'(t_0) \ = \ 
	-s_2h(t_0)+ \p_t\alp'(t_0),
\]
where in the last equality we used $\alp(t_0)=0$. We conclude that
\[
	h'(t_0)+s_2h(t_0)\ = \ \p_t\alp'(t_0)
\]
and 
\[
	 	\n_{d/dt}h(t)\ := \ P_t\,\big(\,h'(t)+s_2h(t)\,\big) \ =\ 0.
\]
\eprf

\subsection{Proof of \refp{Phiinjective}}\label{SS:prPhiinjective}
By \refl{h(t)} the map $\oPhi_{t_1,t_2}$ is equal to the monodromy map of the connection $\n$. Hence, it is an isomorphism. \hfill$\square$ 
\subsection{Proof of \reft{independenceofs}}\label{SS:prindependenceofs}
By \refp{Phiinjective} the maps 
\[
	\oPhi_{s_1,s_1+s_2}:\, \bar\calH_{s_1}^{0,\b} \ \longrightarrow \ \bar\calH_{s_1+s_2}^{0,\b} 
\]
and 
\[
	\oPhi_{s_2,s_1+s_2}:\, \bar\calH_{s_2}^{0,\b} \ \longrightarrow \ \bar\calH_{s_1+s_2}^{0,\b} 
\]
are isomorphisms.  Hence, the map
\[
	\oPhi_{s_1,s_2}\ := \  \oPhi_{s_2,s_1+s_2}^{-1}\circ \oPhi_{s_1,s_1+s_2}:\, 
		\bar\calH_{s_1}^{0,\b} \ \longrightarrow \ \bar\calH_{s_2}^{0,\b} 
\]
is an isomorphism. 

Since $s_1$ and $s_2$ are $V$-generic, $ \bar\calH_{s_i}^{0,\b}$ is canonically isomorphic to $H^{0,\b}_{s_i,V}(M,E)$ ($i=1,2$). We now set
\begin{equation}\label{Phis1s2 definition}
	\Phi_{s_1,s_2}^V\ := \ \oPhi_{s_1,s_2}.
\end{equation}

It  remains to show that $\Phi_{s_1s_2}^V$  satisfies the cocycle condition \refe{cocyclePhi}. Let $s_1,s_2,s_3\in \calF$ be $V$-generic functions. From \refe{oPhit1t2t3} we conclude that 
\eq{4Phi}
	\oPhi_{s_1+s_2,s_1+s_2+s_3}\circ \oPhi_{s_1,s_1+s_2}\ = \ 
	\oPhi_{s_1+s_3,s_1+s_2+s_3}\circ\oPhi_{s_1,s_1+s_3}.
\end{equation}

Recall that  for each $s\in\calF$ the map $\oPhi_{s,s+s_j}$ $(j=1,2,3)$ is induced by the map 
\[
	e^{-s_j}:\,\Ome^{0,\b}(M,E)\ \to \ \Ome^{0,\b}(M,E),
\]
To simplify the notation we will denote  $\oPhi_{s,s+s_j}$ by $e^{-s_j}$. By permuting the indices 1,2,3 in \refe{4Phi} we conclude that the following diagram commutes:
\bigskip
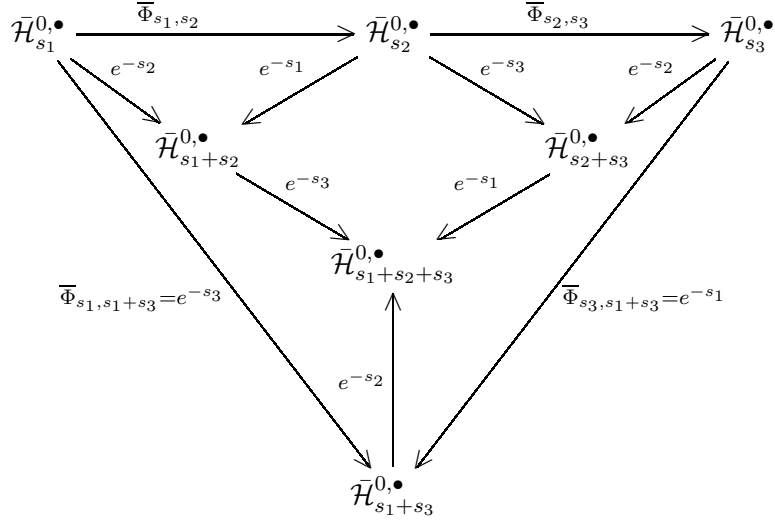
\begin{diagram}
\bar\calH_{s_1}^{0,\b}&\rTo^{\quad \quad \oPhi_{s_1,s_2}\qquad\qquad } &\bar\calH_{s_2}^{0,\b}&\rTo^{\quad \quad \oPhi_{s_2,s_3}\qquad\qquad }&\bar\calH_{s_3}^{0,\b}\\
&\rdTo(2,8)_{\oPhi_{s_1,s_1+s_3}=e^{-s_3}}\rdTo(1,2)^{e^{-s_2}}\ldTo(1,2)^{e^{-s_1}}&
&\rdTo(1,2)^{e^{-s_3}}\ldTo(1,2)^{e^{-s_2}}\ldTo(2,8)_{\oPhi_{s_3,s_1+s_3}=e^{-s_1}}&\\
&\bar\calH_{s_1+s_2}^{0,\b}&&\bar\calH_{s_2+s_3}^{0,\b}&\\
&&\rdTo(1,2)^{e^{-s_3}}\ldTo(1,2)^{e^{-s_1}}&&\\
&&\bar\calH_{s_1+s_2+s_3}^{0,\b}&&\\
&&&&\\
&&\uTo(0,4)^{e^{-s_2}}&&\\
&&&&\\
&&\bar\calH_{s_1+s_3}^{0,\b}&&
\end{diagram}
\medskip
Hence, 
\[
	\oPhi_{s_2,s_3}\circ\oPhi_{s_1,s_2} \ = \ \oPhi_{s_3,s_1+s_3}^{-1}\circ\oPhi_{s_1,s_1+s_3}
	\ = \ \oPhi_{s_1,s_3}.
\]
\hfill$\square$

\section{Proof of \reft{indepofh}}\label{S:prindepofh}

\subsection{The case $h_1^E>h_2^E$}\label{SS:h1>h2}
Suppose first that $h_1^E>h_2^E$. Then $L_2\Ome^{0,\b}(M,E_1)\subseteq L_2\Ome^{0,\b}(M,E_2)$. This inclusion induces  a map of  cohomology 
\[
	i:\, H^{0,\b}_s(M,E_1) \ \longrightarrow \ H^{0,\b}_s(M,E_2).
\]
Since $\tilmu:M\to \grg$ is proper there exists a compact set $K_1\subset M$ and a smooth function 
\[
	\kap{:}\,[0,\infty)\ \to\  [0,\infty)
\]
such that 
\eq{kaph1<h2}
	e^{-2\kap(|\tilmu|^2/2)} \, h^E_1\ <\ h^E_2
\end{equation}
on $M\backslash{}K_1$. By \refl{admissible2}  there exists $\tau\in \calF(M,g^{M},E,h^{E}_{1}+h^{E}_{2})$ with  $\tau\ge\kap$. It follows from \refr{h1<h2} that $\tau\in \calF(M,g^{M},E,h^{E}_{i})$ for $i=1,2$.

Assume now that $s$ is $V$-generic for both bundles $E_1$ and $E_2$. Then it follows from \refp{family of projections} that there exists $t>1$ such that $s+t\tau$ is also $V$-generic.  From \refe{kaph1<h2} we conclude that for any $\alp\in L_2\Ome^{0,\b}(M,E_2)$ 
\[
	e^{-t\tau(|\tilmu|^2/2)}\alp\ \in\ L_2\Ome^{0,\b}(M,E_1).
\]
Since
\[
	\p_{s+t\tau}\circ e^{-t\tau(|\tilmu|^2/2)} \ = \ e^{-t\tau(|\tilmu|^2/2)}\circ \p_{s}
\]
the multiplication by $e^{-t\tau(|\tilmu|^2/2)}$ defines a map
\[
	e^{-t\tau(|\tilmu|^2/2)}:\, H^{0,\b}_{s}(M,E_{2}) \ \longrightarrow H^{0,\b}_{s+t\tau}(M,E_{1}).
\]
By \reft{independenceofs} the map $\alp\mapsto e^{-t\tau}\alp$ induces isomorphisms
\[
	\Phi_{s,s+t\tau} = e^{-t\tau}:\,H^{0,\b}_s(M,E_i) \ \longrightarrow \ H^{0,\b}_{s+t\tau}(M,E_i), \qquad i=1,2.
\]
Since $s$ and $s+t\tau$ are $V$-generic for the bundles $E_1$ and $E_2$, we get an isomorphism 
\[
	\Phi_{s,s+t\tau} = e^{-t\tau}:\,H^{0,\b}_{\bg,V}(M,E_i) \ \longrightarrow \ H^{0,\b}_{\bg,V}(M,E_i), \qquad i=1,2.
\]
Thus we obtain a commutative diagram
\begin{diagram}[height=1.3cm]
H_{\bg,V}^{0,\b}(M,E_1) &\rTo^{\hskip2cm i\hskip2cm} &H_{\bg,V}^{0,\b}(M,E_2)\\
\dTo^{ e^{-t\tau}} & \ldTo^{e^{-t\tau}}&\dTo_{e^{-t\tau}}\\
 H^{0,\b}_{\bg,V}(M,E_1)&\rTo^{i} & H^{0,\b}_{\bg,V}(M,E_2)
\end{diagram}
Since the vertical arrows in this diagram are isomorphisms, it follows that horizontal arrows are also isomorphisms.

 We now set
\[
	\Tet_{h_1^Eh_2^E} ^V\ := \ i:\, H^{0,\b}_{\bg,V}(M,E_1) \ \longrightarrow \ H^{0,\b}_{\bg,V}(M,E_2).
\]

\subsection{Proof of \reft{indepofh} in the general case}\label{SS:prindepofh}
Set $h^E=h_1^E+h_2^E$. Then $h^E>h_i^E$ $(i=1,2)$. By \reft{independenceofs} it is enough to prove  \reft{indepofh} for the function 
\[
	s\ \in \  \calF(M,g^M,E,h^E) \subset \ \calF(M,g^M,E,h_1^E)\,\cap\,\calF(M,g^M,E,h_2^E).
\]
For such $s$ the map
\[
	\Tet_{h_1^Eh_2^E}^V \ :=\ \big(\Tet_{h^E_2h^E}^V\big)^{-1} \circ \Tet_{h_1^Eh^E}^V 
\]
is an isomorphism.  The proof that it satisfies the cocycle condition \refe{cocycleTet} is a verbatim repetition of the arguments in \refss{prindependenceofs} .\hfill$\square$

\section{The Kodaira-type vanishing theorem}\label{S:vanishing}

Let $L$ be a holomorphic line bundle over $M$ and let $\n^L$ be a holomorphic connection on $L$. Recall that $L$ is called positive if the curvature $F^L=(\n^L)^2$ of $\n^L$ is positive, i.e. 
\eq{positive}
	F^L(v,\ov) \ > \ 0, \qquad \text{for all}\quad v\in T^{1,0}M.
\end{equation}
The purpose of this section is to prove the following analogue of the Kodaira vanishing theorem:

\th{kodaira}
Let $L$ be a $G$-equivariant positive  line bundle over a tamed \ka $G$-manifold $M$. For any holomorphic  $G$-equivariant vector bundle $E$ over $M$ and any irreducible representation $V\in \Irr G$, there exists an integer $k_0>0$, such that for all $k\ge k_0$ the $V$-component of the background cohomology 
\[
	H^{0,p}_{\bg,V}(M,E\otimes L^{\otimes k}) \ = \ 0,
\]
for all $p>0$. 
\eth
The rest of this section is occupied with the proof of \reft{kodaira}. First we explain the main idea of the proof.

\subsection{The plan of the proof}\label{SS:planprkodaira}
In \refl{admissibleEk} we show that there exists a function $s$ which is admissible for $E_k=E\otimes{}L^{\otimes k}$ for all $k\in \NN$.  Let $D_k:\Ome^{0,\b}(M,E_k)\to \Ome^{0,\b}(M,E_k)$ denote the Dolbeault-Dirac operator and let $D_{s,k}$ denote its deformation. A direct computation shows that (cf. \refl{D2})
\[
	D_{k,s}^2\ = \ D^2_k\ + \ T_{fv,k} \ - 2i\,\n^{\E_k}_{fv},
\]
where $T_{fv,k}$ is a zero order operator and $\n^{\E_k}$ is the connection on  $\E_k= E_k\otimes \Lam^\b(T^{0,1}M)^*$. The operator $\n^{\E_k}_{fv}$ is a first order differential operator, but for each $V\in \Irr{}G$ its restriction to the $V$-component $\Ome^{0,\b}(M,E_k)^V$ of $\Ome^{0,\b}(M,E_k)$ is a zero order operator. It follows that there exists a function $r_V:M\to \RR$ such that
\eq{Du2a0}
     D_{s,k}^2|_{\Ome^{0,\b}(M,E_k)^V} \ \ge \ D_k^2|_{\Ome^{0,\b}(M,E_k)^V} \ + r_V.
\end{equation}
In \refl{Dsk>} we show by a direct computation that there exists a compact set $K_1\subset M$ such that for all $x\not\in K_1$, we have $r_V(x)> \|B(x)\|+1$, where $B(x)$ is the bundle map defined in the Bochner-Kodaira formula \refe{lechnW}. Then in \refss{prkodaira} we use the Bochner-Kodaira formula and \refe{Du2a0} to show that for large $k\in \NN$
\[
 	D_{s,k}^2|_{\Ome^{0,\b}(M,E_k)^V} \ \ge \ 1, \qquad \text{for all}\quad p>0.
\]
The regularized cohomology $H^{0,p}_{s,V}(M,E\otimes L^{\otimes k})$ is isomorphic to the kernel of $D_{s,k}^2|_{\Ome^{0,p}(M,E_k)^V}$ and, hence, vanish for $p>0$. The theorem follows now from the fact that 
\begin{equation}\label{E:bg<s}
	\dim H^{0,p}_{\bg,V}(M,E\otimes L^{\otimes k}) \le \dim H^{0,p}_{s,V}(M,E\otimes L^{\otimes k}).
\end{equation}

We now present the details of the proof of \reft{kodaira}.

\subsection{The Dolbeault-Dirac operator on $E\otimes L^{\otimes k}$}\label{SS:dirackodaira}
Choose a Hermitian metric $h^E$ on $E$ and a Hermitian metric $h^L$ on $L$. For $k\in \NN$ let $E_k=E\otimes{}L^{\otimes k}$ and let $h^{E_k}$ denote the Hermitian metric on $E_k$ induced by $h^E$ and $h^L$.  Let $\n^E$ and $\n^L$ denote the holomorphic Hermitian connections on $E$ and $L$ and let $\n^{E_k}$ be the induced connection on $E_k$.  We denote by 
\[
	D_k\ := \ \sqrt2\,\big(\,\p+\p^*\big):\,\Ome^{0,\b}(M,E_k)\ \to \Ome^{0,\b}(M,E_k)
\]
the corresponding Dolbeault-Dirac operator on anti-holomorphic differential forms with values in $E_k$. Set
\[
	\E_k\ := \ E_k\otimes \Lam^\b(T\ah M)^*.
\]
The space $\Ome^{0,\b}(M,E_k)$ of anti-holomorphic forms with values in $E_k$ is isomorphic to the space $\Gam(M,\E_k)$ of smooth  sections of $\E_k$.  Let  
\[
	\n^{\E_k} \ := \ \n^{E_k}\otimes1\ + \ 1\otimes \n^{LC}
\] 
denote the connection on $\E_k$ induces by $\n^{E_k}$ and the Levi-Civita connection $\n^{LC}$ on $ \Lam^\b(T\ah M)^*$.

\subsection{The action of $F^L$}\label{SS:FLaction}
Consider an endomorphism $\lam(F^L)$ of $\calE_k$ defined by the formular (cf. \cite[\S3]{BeGeVe})
\begin{equation}\label{E:bfc2}
        \lam(F^L)\,\alp \ = \ \sum_{i,j}\, F^L(w_i,\ow_j)\, w_i\wedge \iot_{\ow_j}(\alp), \qquad  \alp\in \calE_k,
\end{equation}
where $\{w_1\nek w_n\}$ is an orthonormal basis of $T^{1,0}M$. 

The positivity assumption \refe{positive} immediately implies  (cf. for example  formula (1.5.19) of \cite{MaMarinescu_book}) the following
\lem{bfc>0}
For all $p>0$ the restriction of $\lam(F^L)$ to the space $\Ome^{0,p}(M,E_k)$ is a strictly positive operator
\[
	\lam(F^L)\big|_{\Ome^{0,p}(M,E_k)}\ > \ 0, \qquad p>0.
\]
\elem

\subsection{The Bochner-Kodaira formula}\label{SS:Lichnerowicz}
It follows from the Bochner-Kodaira formula \cite[Theorem~3.71]{BeGeVe}, that
\begin{equation}\label{E:lechnW}
        D_k^2 \ = \ \Del_k^{0,1} \ +\, k\,\lam(F^L) \ + \ B,
\end{equation}
where 
\[
	\Del_k^{0,\b} \ = \ \big(\n^{0,1}\big)^*\n^{0,1}
\]
is the generalized Laplacian on $\E_k$ and $B\in\End(\E_k)$ is independent of $k$.

\subsection{Admissible functions for  $E_k$}\label{SS:admissiblek}
For $\u\in\grg$ let $u$ denote the corresponding vector field on $M$, cf. \refe{v}.
We denote by $\L^\E_\u$, $\L^{\E_k}_\u$, and $\L^L_\u$ the infinitesimal action of $\u$  on $\Gam(M,\E)$,  $\Gam(M,\E_k)$, and  $\Gam(M,L)$ respectively. As in \refss{rescaling},  we set
\eq{muEL}
	 \mu^{\E_k}(\u) \ := \ \n^{\E_k}_{u}-\L^{\E_k}_\u, \quad  \mu^L(\u) \ := \ \n^L_{u}-\L^L_\u.
\end{equation}
Then
\eq{muEn=}
	\mu^{\E_k}(\u) \ = \ \mu^\E(\u) \ + \ k\,\mu^L(\u).
\end{equation}
By Kostant formula (\cite{Kostant70}, \cite[(1.13)]{TianZh98}) the moment map $\mu$ is related to $\mu^L$ by 
\eq{muL}
	\mu^L(\u) \ = \ {2\pi}i\,\<\mu,\u\>.
\end{equation}
Combining \refe{muEn=} and \refe{muL} we obtain
\[
    \mu^{\E_k}(\u) \ = \ \mu^\E(\u) \ + \ {2\pi k}i\,\<\mu,\u\>.
\]
In particular, when $\u=\v$ is given by \refe{bfv(x)}, we obtain 
\eq{muEN=muv}
    \mu^{\E_k}(\v) \ = \ \mu^\E(\v) \ + \ {2\pi k}i\,|\v|^2.
\end{equation}
Hence, from \refe{muEL} we get
\eq{nEk=}
 	\n^{\E_k}_\v\ = \ \L^{\E_k}_\v\ +\ \mu^\E(\v)\ + \ {2\pi k}i\,|\v|^2.
\end{equation}

\lem{admissibleEk}
There exists a function $s:[0,\infty)\to [0,\infty)$ with the following properties:
\begin{enumerate}
\item it  is admissible for $E_k$ for all $k\in \NN$;
\item there exits a compact set $K\subset M$ such that for all $x\in M\backslash{}K$ 
\eq{f>B}
	f(x)\ := \ s'\big(|\mu(x)|^2/2\big) \ > \ \frac{\sqrt{2\|B(x)\|+2}}{|v(x)|}.
\end{equation}
Here $B\in \End(\E)\subset \End(\E_k)$ is defined in \refe{lechnW}.
\end{enumerate}
\elem
\prf
Fix a function $\tau:M\to [0,\infty)$ with $\lim_{M\ni x\to \infty}\tau(x)=\infty$. By  \refl{admissible2} we can choose a function $s$ admissible for $(E,h^E)$ such that for all $x\in M\backslash{K}$
\eq{f>}
	f(x)\ := \ s'\big(|\mu(x)|^2/2\big) \ > \ \max\Big\{\,\frac{|\v(x)|^2}{|v(x)|^2}\,\tau(x),\,\frac{\sqrt{2\|B(x)\|+2}}{|v(x)|}\,\Big\}.
\end{equation}
Let $\nu$ be as in \refe{nu} and let 
\[
    \nu_k=|\v|+\|\nLC v\|+\|\mu^{\E_k}(\v)\|+|v|+1.
\]
Then from  \refe{muEN=muv} and \refe{f>} we obtain
\[
	\nu_k \ \le \nu \ +\  2\pi k |\v|^2 \ \le \ \nu \ + \ \frac{2\pi k f|v|^2}\tau.
\]
Hence, 
\eq{admissiblek}
	        \frac{f^2|v|^2}{    |df||v|+f\nu_k+1 }  \  \ge \ \frac{f^2|v|^2}{    |df||v|+f\nu+1+\frac{2\pi k f^2|v|^2}\tau }
\end{equation}
Using \refe{admissible} and the fact that $\lim_{M\ni x\to \infty}\tau(x)=\infty$, we conclude from \refe{admissiblek} that 
\[
	\lim_{M\ni x\to \infty}\, \frac{f^2|v|^2}{    |df||v|+f\nu_k+1 }  \  = \ \infty. 
\]
\eprf

\subsection{Deformation of $D_k$}\label{SS:defDk}
Let $s$ be a function admissible for $E_k$ for all $k\in \NN$ and let $\phi(x)=s(|\mu|^2/2)$. Let 
\[
	D_{s,k}\ := \ \sqrt2\,\big(\,\p_s+\p_s^*\big):\,\Ome^{0,\b}(M,E_k)\ \to \Ome^{0,\b}(M,E_k),
\]
where as in \refe{deformedDobeault}
\[
	\p_s \alp\ = \ e^{-\phi}\circ\p\circ e^\phi\,\alp 
	\ = \ \p\alp\ + \ f\,\p\big(\,|\mu|^2/2\,\big)\wedge\alp, \qquad \alp\in \Ome^{0,\b}(M,E_k).
\]

Consider the operator
\eq{Av}
        A_u \ = \ \sum c(e_i)\, c(\nLC_{e_i}u):\, \E_k \ \to \ \E_k.
\end{equation}
\lem{D2} 
The following equality holds 
\eq{D2} 
    D_{s,k}^2 \ = \ D_k^2 \ + \ f^2|v|^2 \ + \ {\i}A_{fv} \ - \ 2{\i}\n^{\E_k}_{fv}.
\end{equation}
\elem
The proof of the lemma is a straightforward calculation, cf. \cite[Lemma~9.2]{Br-index} or \cite[Theorem~1.6]{TianZh98}.

\lem{Dsk>}
Let $s$ be the admissible function constructed in \refl{admissibleEk}. There exists a smooth function $r_V:M\to \RR$ and a compact set $K_1\subset M$ such that 
\eq{rV>B}
	r_V(x)\ >\ \|B(x)\|+1.
\end{equation}
for all  $x\in M\backslash{}K_1$ and  
\eq{Du2a}
     D_{s,k}^2|_{\Ome^{0,\b}(M,E_k)^V} \ \ge \ D_k^2|_{\Ome^{0,\b}(M,E_k)^V} \ + r_V.
\end{equation}
\elem
\prf
Since $\|c(v)\|=|v|$ and $\|c(e_i)\|=1$, we
have
\eq{Av2}
    \|A_{fv}\| \ \le \ \sum_i\, \|\nLC_{e_i}(fv)\| \ \le \            C\, \Big(\, |df|\, |v| + f\, \|\nLC v\|\, \Big),
\end{equation}
for some constant $C>0$.

For $\u\in \grg$ let $\calL_\u^V:V\to V$ denote the action of $\u$ on $V$. Then there exists a constant $c_V$ such that for any $G$-invariant scalar product on $V$ we have
\[
	\|\calL_\u\|\ \le \ c_V|\u|.
\]
Consider the space 
\eq{GamV}
	\Ome^{0,\b}(M,E_k)^V\ = \ \Gam(M,\E_k)^V \ := \ \Hom_G\big(\,\Gam(M,\E_k),V\,\big)\otimes V.
\end{equation}
The restriction of the operator $\calL^{\E_k}_\u$ to $\Gam(M,\E_k)^V$ decomposes as $\calL^{\E_k}_\u=1\otimes\calL^V_\u$. Hence 
\eq{LEk<}
    \big\|\, \L^{\E_k}_\u|_{\Gam(M,\E_k)^V}\, \big\| \ \le \ c_V |\u|.
\end{equation}
Combining, \refe{D2}, \refe{Av2}, \refe{LEk<} and \refe{nEk=}, we obtain
\begin{multline}\label{E:Dsk2>}
    D_{s,k}^2|_{\Ome^{0,\b}(M,E_k)^V} \ \ge \
    D_k^2|_{\Ome^{0,\b}(M,E_k)^V}  \\ + \  f^2\, |v|^2 \ - \ \lam_V\, \Big(\, |df|\, |v|  \ + \ f\, \big(\, |\v|+\|\mu^\E(\v)\|+\|\nLC v\|\, \big) \, \Big)\ + \ 4\pi kf\,|\v|^2,
\end{multline}
where $\lam_V=\max\{1,c_V,C\}$. 
Set 
\[
	r_V\ := \ f^2\, |v|^2 \ - \ \lam_V\, \Big(\, |df|\, |v|  \ + \ f\, \big(\, |\v|+\|\mu^\E(\v)\|+\|\nLC v\|\, \big) \, \Big)
	\ + \ 4\pi kf|\v|^2.
\]
Then \refe{Dsk2>} is equivalent to \refe{Du2a}. From \refe{admissible} we conclude that there exists a compact set $K_1\subset M$ such that 
\eq{rV>}
	 r_V(x)\ \ge \frac12f^2(x)\, |v(x)|^2 
\end{equation}
for $x\not\in K_1$. 

Let $K$ be the compact set defined in \refl{admissibleEk}. We can assume that $K_1\supset K$. Then the inequality \refe{rV>B} follows  form  \refe{f>B}.
\eprf

\subsection{Proof of \reft{kodaira}}\label{SS:prkodaira}

From  \refl{bfc>0} we see that for large enough integer $k$, for every $x\in K_1$, and  for every $p>0$ 
\[
	k\,\lam(F^L)(x) \ + \ B(x)\ + \ r_V(x) \ \ge \ 1.
\] 
From Lemmas~\ref{L:bfc>0} and \ref{L:Dsk>} we conclude that for any $k>0$ and $x\not\in K_1$ we have
\[
	k\,\lam(F^L)(x) \ + \ B(x)\ + \ r_V(x) \ \ge\ B(x)\ +\ r_V(x)\ \ge \ 1.
\]
Combining these two inequalities with  the Bochner-Kodaira formula \refe{lechnW} and \refl{Dsk>} we conclude that for $p>0$
\[
	D_{s,k}^2|_{\Ome^{0,\b}(M,E_k)^V} \ \ge \ \Del_k^{0,\b}|_{\Ome^{0,p}(M,E_k)^V} \ + \ 1 \ \ge \ 1.
\]
Hence, 
\eq{KerD=0}
	\Ker{}D_{s,k}^2|_{\Ome^{0,p}(M,E_k)^V}\ =\ 0, \qquad \text{for all} \ \ p>0.
\end{equation}

The space $H^{0,p}_{\bg,V}(M,E_k)$ is equal to the cohomology of the complex  $\big(\Ome^{0,\b}(M,E_k)^V,\p_s\big)$. Thus $H^{0,p}_{\bg,V}(M,E_k)$ is isomorphic to the kernel of the restriction of $D_{s,k}^2$ to $\Ome^{0,p}(M,E_k)^V$. \reft{kodaira} follows now from \refe{KerD=0} and \refe{bg<s}. 
\hfill$\square$

\appendix
\section{Dependence of the spectrum of $D_s$ on $s$}\label{S:holomorphicB}

In this appendix we study the spectral properties of the family of operators $D^V(t)$ introduced in \refss{backgroundharmonic} and prove \refp{family of projections}. Throughout the appendix we use the notation introduced in \refs{prindependenceofs}.

\subsection{Small values of $t$}\label{SS:smallt}
We start with showing that the admissible function $s(t)$ are $V$-generic for small $t$. 

\lem{smallt}
There exists $\del>0$ such that for all $t\in [0,\del)$ we have 
\begin{equation}\label{E:smallt}
	\dim \Ker D^V(t) \ = \ \dim   H^{0,\b}_{\bg,V}(M,E).
\end{equation}
Moreover, $P_t$ ($t\in [0,\del)$) is a continuous family of projections. 
\elem
\prf
It is shown in the proof of Proposition~10.5 of \cite{Br-index} that the spectrum of $D^V(t)$ is discrete. Hence, there exists an $\eps>0$ such that 
\begin{equation}\label{E:D>eps}
	D^V(0)^2\big|_{\IM(\Id-P_0)} \ >\ \eps.
\end{equation}

The same argument as in the proof of \refl{Dsk>} shows that there exists a compact set $K\subset M$ and a smooth function $r_V(x)$ such that $r_V(x)>1$ for all $x\in M\backslash{K}$ and 
\[
	D^V(t)^2 \ \ge \ D^V(0)^2\ + t\,r_V(x).
\]
Thus there exists $m>0$ such that $r_V(x)\ge -m$ for all $x\in M$. Then 
\eq{D>D-tm}
	D^V(t)^2\ \ge\  D^V(0)^2-tm.
\end{equation}
Set
\[
	\del\ := \ \frac{\eps}{2m}.
\]
Combining \refe{D>eps} with \refe{D>D-tm} we obtain
\eq{D>0}
	D^V(t)^2\big|_{\IM(\Id-P_0)} \ >\ \eps/2, \qquad\text{for all}\quad t\in [0,\del)
\end{equation}
Hence, 
\[
	\dim \Ker D^V(t)\ = \ \dim \Ker D^V(t)^2 \ \le \ \dim \Ker D^V(0), \qquad t\in [0,\del)
\]
Since $s$ is a $V$-generic function 
\[
	\dim\Ker D^V(0)\ = \ \dim H^{0,\b}_{\bg,V}\ \le \ \Ker D^V(t)
\] 
for all $t$. Hence, we obtain \refe{smallt}.

Let $\gam(\tau)= \frac\eps3e^{i\tau}$. Then there are no non-zero eigenvalues of $D^V(t)^2$ on or inside $\gam$ for all $t\in[0,\del)$. Hence, 
\[
	P_t\ = \ \frac1{2\pi i}\int_\gam\big(\lam-D^V(t)^2\big)^{-1}\,d\lam
\]
depends continuously on $t$. 
\eprf

\subsection{Continuity of $P_t$ for $t>0$}\label{SS:Ptt>0}
We now study the dependence of the spectrum of $D^V(t)$ on $t$ for $t>0$. First, we show that $D^V(t)^2$ is a holomorphic family of operators. For this we represent 
\[
	D^V(t)^2\ = \ D^V(t_0)^2 \ + \ (t-t_0)\,S_1\ +\ (t-t_0)^2\,S_2
\]
and apply Theorem~4.12 of \cite[Ch.~VII]{Kato}. We now make this argument more precise.

\lem{typeB}
The family of operators $D_V(t)^2,\ t>0$ is a holomorphic family of type $(B)$ in the sense of Kato \cite[Ch.~VII\, \S4.4]{Kato}.
\elem
\prf
Fix $t_0>0$ and let 
\[
	f_i\ := \ s_i'\big(|\tilmu|^2/2), \qquad i=1,2
\]
As in \refl{Dsk>} one shows that there exists a smooth function $r_V:M\to \RR$ such that 
\eq{DVt02}
	D^V(t_0)^2\ = \ D^V(0)^2\ +\ r_V(x)
\end{equation}
and 
\begin{equation}\label{E:rV>12f2v}
	r_V(x)\ > \ \frac12\, f_2(x)^2|v(x)|^2
\end{equation}
for all $x$ outside of a compact subset of $M$ (the last inequality is obtained the same way as \eqref{E:rV>}). 

Set 
\[
	a_1\ := \  \max_{x\in M} \left(\,\frac12\, f_2(x)^2|v(x)|^2-r_V(x)\,\right).
\]
Then it follows from \refe{DVt02} that for every form $\ome$ in the domain of $D^V(0)^2$ we have
\eq{DVome>}
	\big\<D^V(t_0)^2\,\ome,\ome\,\big\> \ > \  \left(\,\frac12\, f_2(x)^2|v(x)|^2-a_1\,\right)\, \|\ome\|^2.
\end{equation}

As in \refl{D2} one shows that 
\eq{DVt2}
	D^V(t)^2\ = \  D^V(t_0)^2 \ + \ (t-t_0)^2f_2^2|v|^2 \ + \ {\i}(t-t_0)A_{f_2v} \ - \ 2{\i}(t-t_0)\n^{\E}_{f_2v},
\end{equation}
where $\E=E\otimes \Lam^\b(T\ah M)^*$.  Set 
\[
	S_1\ := \ A_{f_2v}-2\i \n^{\E}_{f_2v}\big|_{\Ome^{0,\b}_V(M,E)}, \quad S_2\ := \ f_2^2|v|^2.
\]
Then 
\begin{equation}\label{E:D2=D+S1+S2}
	D^V(t)^2\ = \  D^V(t_0)^2 \ + \ (t-t_0)\,S_1\ +\ (t-t_0)^2\,S_2.
\end{equation}

From \refe{DVome>} we obtain
\begin{equation}\label{E:S2<D}
	\big|\, \<\,S_2\ome,\ome\>\,\big| \ \le \ 2\, \big\<D^V(t_0)^2\,\ome,\ome\,\big\> + 2a_1\|\ome\|^2,
\end{equation}
for every $\ome$ in the domain of $D^V(t_0)^2$.

It is shown in the proof of \refl{Dsk>} that $S_1$ is a bundle map and that there exists a compact set $K_1\subset M$ and a constant $\lam_V$ such that 
\[
	\big\|\, S_1(x)\big\| \ \le \  \lam_V\, \Big(\, |df_2|\, |v|  \ + \ f_2\, \big(\, |\v|+\|\mu^\E(\v)\|+\|\nLC v\|\, \big) \, \Big)
\]
for all $x\not\in K_1$.  From \refe{admissible} we now conclude that there exists a compact set $K\supset K_1$ such that 
\begin{equation}\label{E:S1<}\notag
	\big\|\, S_1(x)\big\| \ \le \ \frac12f_2^2(x)|v(x)|^2
\end{equation}
for $x\not\in K$. Set 
\[
	a_2\ := \ \max_{x\in K} \|S_1(x)\|, 
\]
then 
\begin{equation}\label{E:S1<2}\notag
	\big\|\, S_1(x)\big\| \ \le \ \frac12f_2^2(x)|v(x)|^2 \ +\ a_2. 
\end{equation}
Comparing with \refe{DVome>} and setting $a=a_1+a_2$ we obtain
\begin{equation}\label{E:S1<D}
	\big|\,\< S_1(x)\ome,\ome\>\,\big| \ \le \  a\|\ome\|^2\ +\ \big\<D^V(t_0)^2\,\ome,\ome\,\big\> , 
\end{equation}
for all $\ome$ in the domain of $D^V(t_0)^2$. 

Combining \refe{D2=D+S1+S2} with \refe{S2<D} and \refe{S1<2} and using Theorem~4.12 of \cite[Ch.~VII]{Kato}, we conclude that there exist an $\eps>0$ such that $D^V(t)^2$ is a holomorphic family of operators of type (B) on the interval $(t_0-\eps,t_0+\eps)$. Since this result holds for every $t>0$ the lemma is proven. 

\subsection{Proof of \refp{family of projections}}\label{SS:prfamily of projections} 
Since $D^V(t)^2$ is a holomorphic family of type (B) for $t>0$, the eigenvalues of $D^V(t)^2$  depend analytically on $t>0$ by the result of Section~VII.4.6 of \cite{Kato}. It follows that there are finitely many eigenvalues of $D^V(t)^2$ which a identically equal to zero and a discrete sequence of positive numbers $t_1, t_2, ...$ such that  the rest of the eigenvalues do not vanish for $t\not\in \{t_1,t_2,...\}$. Moreover, the eigenfunctions corresponding to zero eigenvalue also can be chosen to be analytic functions of $t$. Let $P_t$ denote the orthogonal projection onto the span of these eigenfunctions. Then $P_t$ ($t>0$) is a holomorphic family of projections and, in particular, depends continuously on $t$. Combination of  this result with \refl{smallt} proves \refp{family of projections}.
\eprf

\def\cprime{$'$} \newcommand{\noop}[1]{} \def\cprime{$'$}
\providecommand{\bysame}{\leavevmode\hbox to3em{\hrulefill}\thinspace}
\providecommand{\MR}{\relax\ifhmode\unskip\space\fi MR }
\providecommand{\MRhref}[2]{%
  \href{http://www.ams.org/mathscinet-getitem?mr=#1}{#2}
}
\providecommand{\href}[2]{#2}

\end{document}